# THE EMERGENCE OF THE DETERMINISTIC HODGKIN–HUXLEY EQUATIONS AS A LIMIT FROM THE UNDERLYING STOCHASTIC ION-CHANNEL MECHANISM

By Tim D. Austin

*University of California at Los Angeles*

In this paper we consider the classical differential equations of Hodgkin and Huxley and a natural refinement of them to include a layer of stochastic behavior, modeled by a large number of finite-state-space Markov processes coupled to a simple modification of the original Hodgkin–Huxley PDE. We first prove existence, uniqueness and some regularity for the stochastic process, and then show that in a suitable limit as the number of stochastic components of the stochastic model increases and their individual contributions decrease, the process that they determine converges to the trajectory predicted by the deterministic PDE, uniformly up to finite time horizons in probability. In a sense, this verifies the consistency of the deterministic and stochastic processes.

## 1. Introduction: Ion channels of excitable membranes.

Most neurons in most organisms have an *axon*: a long, narrow conduit connecting the central, roughly spherical part of the cell (the *soma*) to a network of smaller branches and ultimately to the *synapses*, which form connections with other neurons (principally at branched projections from the latter called *dendrites*). The axon connects the soma to synapses that may be a great distance away (often several cm) relative to the size of the soma or the diameter of the axon (typically a few $\mu$m). The function of the neuron relies partly on its ability to transmit signals from the soma to other neurons over this long distance via the axon.

The nature of these signals had begun to become clear during the 1930s, but only with Hodgkin and Huxley's (Nobel Prize-winning) work on the mechanism of signal transmission in the squid giant axon in the early 1950s were the first foundations laid of an accurate mathematical model of their behavior (see [12]).

Since then Hodgkin and Huxley's original analysis has been extended and refined repeatedly. The mathematics underlying the resulting models has been studied for the sake of both more accurate numerical modeling and better theoretical understanding. In particular, Hodgkin and Huxley's empirical, deterministic model has been refined to a model of the axon in which the relevant behavior arises from the combined contributions of a large number of small stochastic components.

Since the present paper is primarily about mathematics, we will assume familiarity with the basic physiological origins of the deterministic and stochastic models (in particular, the working of voltage-dependent ion channels and the action potential). In these terms we will give a brief motivation for the two models in Section 2.1. However, once we have reached the definitions of the equations themselves further references to the physiology will be peripheral, and not important for understanding the paper. A thorough treatment of this physiology can be found in Hille's classic text [11], while a more mathematical description of various such models can be found in Cronin [2].

It is worth noting that, while the deterministic mathematical model has been studied intensively (particularly for numerical computation purposes), the results for the stochastic model are fairly few. As far as I am aware, the pure mathematics behind the stochastic model considered below has never been worked out in detail. In this paper we will prove an existence theorem for that stochastic process, and, more interestingly, the convergence of its various components (such as the function giving the membrane potential along the axon at a particular instant in time) to their counterpart trajectories in the deterministic theory, uniformly up to finite time horizons in probability.

(Results analogous to this have been obtained by Fox and Lu [9] (building on a simulation method of DeFelice and Isaac [3]) for the case in which the membrane potential is assumed constant along the entire length of the axon at each instant. In this case the partial differential equations we will encounter simplify to ordinary differential equations, coupled with a finite number of discrete stochastic processes that can then be studied using the standard methods of Fokker–Planck and Langevin equations. In fact, this simpler case corresponds more closely to the original experimental set up of Hodgkin and Huxley, in which a fine conducting silver wire was inserted along the axon, causing the membrane potential to adjust to a single common value along the axon effectively instantaneously.)

The consequences of the general stochastic model have received increased interest in recent efforts, first by Chow and White [1] and then Faisal, White and Laughlin [8], to estimate how much noise the actual stochastic nature can introduce to a real neuron (behavior that would not appear in the deterministic approximation) and what constraints this places on the size of the axon if it is to function reliably. We will remark more on this briefly in Section 5.

REMARK. When a suitable stimulus is applied (e.g., from the soma at one end of the axon), exceeding a certain threshold, the trajectory of the potential difference along the axon evolves through a family of subthreshold configurations into an action potential. After moving away from its point of origin, this trajectory asymptotically takes the special form of a traveling wave. Although the possibility of such a traveling wave is key to the axon's ability to transmit a signal, we will not refer to it again in this paper. Our subsequent convergence results require only the existence of *some* sufficiently regular time evolution of the system given appropriate initial conditions.

**2. The mathematical models.** In this section we describe the precise mathematical models that we will study. We will assume standard notions from stochastic analysis and PDE.

2.1. *Basic components.* This subsection assumes some knowledge of the physiology of axons; the disinterested reader may skip to the definitions of the equations in the next subsection without impediment.

In microscopic detail, the instantaneous electrical state of the axon depends on the locations of all the ions in solution inside and outside the axon, on the locations and internal states of any molecular mechanisms at work in the axon (the ion channels, in particular), and on various other components of the system. As usual, we do not actually work at this level of detail, but instead make a number of simplifications. However, there is some choice in this procedure. We will see that heuristically the two different models to be studied arise from two different such approximations, one coarser than the other: in particular, the stochastic model describes the working of individual ion channels, whereas the deterministic model "averages out" their behavior, involving instead functions that describe the proportion of those channels in a small neighborhood of a point that are in each possible state.

Before explaining this difference in more detail, we will describe some simple approximations that are made in both cases. First, although the axon is described above as a tube (and actually has a membrane with considerable molecular structure of its own and further cellular components within the axonal fluid), its diameter is so small compared with its length (typically less than 10 $\mu$m compared with a few cm) that all the relevant quantities vary only negligibly across it, and so we simplify the geometric description of the axon to an interval $I = [-\ell, \ell]$. We write $I^\circ$ for the interior $(-\ell, \ell)$ of this interval.

It is worth noting that voltage-dependent ion channels embedded in a cellular membrane seem to be involved in physiological processes other than the working of axons, sometimes in a pattern spread over a nontrivial area, rather than the approximately-linear distribution we are assuming here. It is possible that the analysis that follows below could be adapted to this case, but I have not tried this. For details of such physiology see Hille [11].

Second, in both models we also approximate the distribution of the individual ions in space by continuous concentrations, to be described by suitable PDE. However, it is possible to go further and avoid altogether the need to work with separate data for each different kind of ion. To do this we deal instead only with the variation of membrane potential along the axon. That this contains all the information we need follows from the further assumption that the flows of the different ions across the axonal membrane, although enough to give rise to the relevant changes in the membrane potential, are negligible compared with the concentration levels that remain both inside and outside the axon. This assumption guarantees that the concentration gradients change only negligibly during the working of the axon, and therefore that the effect of ion influx and efflux on membrane potential is correctly described by an equivalent driving potential for each kind of ion. We will omit the relevant mathematical working to establish this description here; see page 37 of Hille [11].

Thus the state of our system is partly described by a function $v\colon I \to \mathbb{R}$ giving the value of the membrane potential at each point along the axon. Since ions can diffuse along the axon, the variation of this function with time will also exhibit diffusive behavior, allowing us to impose certain regularity conditions on it. For simplicity we will assume that the diffusivity constant is 1 throughout this paper; a simple scaling of $I$ recovers the general case.

It will turn out to be suitable that we assume $v$ Lipschitz and in the Sobolev space $H^1$. In fact, we will assume a little more for technical simplicity. In a real axon the equilibrium potential need not be zero; however, we can and will shift our origin so that we can treat it as zero. We do this because we will find it helpful to impose the condition that our functions vanish at $\pm\ell$ and so to restrict attention to potential difference functions in $H_0^1(I)$. If we do not make this change of origin, then our results for such functions will certainly still be valid; the problem, rather, is that they will no longer apply to the biophysically interesting situation.

We will write $\mathbf{v}_t$ for the potential difference function at time $t$ in the deterministic model, and $\mathbf{V}_t$ for that in the stochastic model. In both cases these should evolve following a continuous trajectory in $H_0^1(I)$ (in the stochastic case, this means as a process in this space with continuous sample paths).

Next we must decide how to model the ion channels and their effect on the membrane potential; it is here that our two models of the axon will diverge.

We do restrict ourselves in both models to the case in which all ion channels are identical, and can be in any of a finite set $E$ of possible channel states. On the other hand, it turns out that there is no great increase in difficulty if we allow several of the possible states $\xi \in E$ to allow the passage of ions with different conductivities. In reality there are different kinds of channel for different kinds of ion, but one finds at once that the resulting mathematical descriptions differ only in notation. Also, in practice there is a constant leak conductivity (corresponding to ion flow across the axonal membrane other than through channels); for simplicity we ignore this also, as we may assume that it has been absorbed by a suitable modification of the conductivities of our channels.

As already mentioned, to each kind of ion there corresponds an equivalent potential difference which "drives" the passage of those ions either into or out of the axon through their corresponding ion channels. Given our treatment of all channels as identical, as described above, in our model these driving potentials $v_\xi$ will actually correspond to the different possible channel states $\xi \in E$ (it will follow at once from the form of the equations that in this arrangement suitable values of the $v_\xi$ can be taken as sums of those corresponding to different types of ion, and that in the case of a state $\xi$ that allows no ions to flow, the value we give to $v_\xi$ will be of no consequence). In our stochastic model a channel at position $x$ will jump between states $\xi, \zeta$ at random at rates $\alpha_{\xi,\zeta}(V)$, where $V$ is the value of the potential difference at the relevant point $x$. We will write $\alpha_\xi$ for the total rate of leaving state $\xi$: $\alpha_\xi = \sum_{\zeta \in E} \alpha_{\xi,\zeta}$.

We will assume that the functions $\alpha_{\xi,\zeta}$ are all smooth and take values between two fixed constants in $(0, \infty)$ (this certainly holds in the actual models that are used). I do not know to what extent this condition could be weakened in what follows; certainly some regularity is needed, and we will use the finiteness of $\mathrm{Lip}(\alpha_{\xi,\zeta})$ explicitly.

Write $v_- = \min_{\xi \in E} v_\xi$, $v_+ = \max_{\xi \in E} v_\xi$, and assume that $v_- < 0 < v_+$ (this is also true in real axons).

Now we can describe our stochastic model; in fact there will be one such model for each $N \in \mathbb{N}$. In the $N$th member of this sequence, the axon is populated by $\lceil 2N\ell \rceil - 2$ channels at positions $\frac{1}{N}(\mathbb{Z} \cap NI^\circ)$, and each has the normalized ion conductivities $\frac{1}{N}c_\xi$ corresponding to the states $\xi \in E$ (the values $c_\xi \geq 0$ being fixed independent of $N$).

We will generally write $\boldsymbol{\Xi}_t$ for the configuration of all the channels in such a stochastic model; this is in the state space $E^{\mathbb{Z} \cap NI^\circ}$. If we want to make $N$ explicit, we include it as a superscript, as in $\boldsymbol{\Xi}_t^{(N)}$ and $\mathbf{V}_t^{(N)}$.

Our deterministic model arises heuristically as the limit of the stochastic model with very many very small ion channels; that is, for large $N$. In the deterministic model we introduce a new family of functions, $p_\xi \in \mathrm{Lip}(I, [0, 1])$ for $\xi \in E$, that replicates the role of the individual-channel configurations $\Xi \in E^{\mathbb{Z} \cap NI^\circ}$. The value $p_\xi(x)$ is to be interpreted as "the proportion of those channels in a small neighborhood of the point $x$ that are in state $\xi$"; we will see that at all times $\sum_{\xi \in E} p_\xi = 1$.

We will write $\mathbf{p}_{\xi,t}$ for the proportion functions at time $t$; these should all evolve following continuous paths in $\mathrm{Lip}(I, [0,1])$.

REMARKS ON NOTATION. Henceforth we will use the Sobolev spaces $H_0^1$ and $H^{-1}$ without further comment. Many good treatments of Sobolev and other function spaces are available in standard texts on PDE; see, for example, Chapter 5 and Section 7.1 of Evans [7].

Given a nonnegative integrable function $f \in L^1(I)$ and writing $\mu$ for Lebesgue measure on $I$, we denote by $\mu \llcorner f$ the indefinite-integral measure:

$$\mu \llcorner f(Y) = \int_Y f \, d\mu.$$

We will also sometimes regard such measures as bounded linear functionals on one or other function space. Given a function $g$ in such a space and a functional $\mu$, we write $\langle g, \mu \rangle$ for the evaluation in the obvious way.

We will write $D$ for differentiation of differentiable functions on $I$ and $\Delta$ for the one-dimensional Laplacian $D^2$, and will use the notation $\chi_E$ for the indicator function of a set $E$.

2.2. *The deterministic equations.* Henceforth suppose that we have $v_0 \in H_0^1$ with $v_- \le v_0 \le v_+$ and a family $(p_{\xi,0})_{\xi \in E}$ of Lipschitz functions $I \to [0,1]$ such that $\sum_{\xi \in E} p_{\xi,0} = 1$ everywhere; these are our initial conditions for the deterministic model described in the previous subsection. We also fix now and for the rest of the paper a finite but arbitrary time horizon $T > 0$.

We are now ready to make the following definition. Note that we are using implicitly a suitable notion of weak solution for our PDE (since we ask only that the time derivative of the trajectory be in the space of functionals $H^{-1}$).

DEFINITION 2.1. A continuous function $\mathbf{v} : [0,T] \to H_0^1(I)$ and a family $(\mathbf{p}_\xi)_{\xi \in E}$ of continuous functions $\mathbf{p}_\xi : [0,T] \to \mathrm{Lip}(I, [0,1])$ will be said to *satisfy the generalized deterministic Hodgkin–Huxley equations* (D) *with initial conditions* $v_0$, $p_{\xi,0}$ if

- (Regularity)

$$\frac{d}{dt}\mathbf{v} \in L^2_{H^{-1}(I)}[0,T],$$

$$\frac{d}{dt}\mathbf{p}_\xi \in L^\infty_{C(I)}[0,T] \qquad \forall \xi \in E;$$

- (Dynamics: PDE)

$$\frac{d}{dt}\mathbf{v}_t = \Delta\mathbf{v}_t + \sum_{\xi\in E} c_\xi \mathbf{p}_{\xi,t}\cdot(v_\xi - \mathbf{v}_t) \qquad \forall t\in[0,T]$$

  [we will refer to this equation as (D-PDE)];
- (Dynamics: proportions)

$$\frac{d}{dt}\mathbf{p}_{\xi,t} = \sum_{\zeta\in E\setminus\{\xi\}} \left((\alpha_{\zeta,\xi}\circ\mathbf{v}_t)\cdot\mathbf{p}_{\zeta,t} - (\alpha_{\xi,\zeta}\circ\mathbf{v}_t)\cdot\mathbf{p}_{\xi,t}\right) \qquad \forall\xi\in E,\ t\in[0,T]$$

  [we will refer to this system of equations as (D-prop)];
- (Initial conditions: PDE)

$$\mathbf{v}_0 = v_0;$$

- (Initial conditions: proportions)

$$\mathbf{p}_{\xi,0} = p_{\xi,0} \qquad \forall\xi\in E;$$

- (Boundary conditions: PDE only)

$$\mathbf{v}_t(\pm\ell) = 0 \qquad \forall t\in[0,T].$$

REMARK. It follows at once by adding the relevant differential equations that the sum $\sum_{\xi\in E}\mathbf{p}_{\xi,t}$ is constant, and so is always equal to 1 everywhere; this means we remain safe in our interpretation of $\mathbf{p}_{\xi,t}$ as the proportion of channels in a particular state.

2.3. *The stochastic equations.* We carry over the PDE initial condition $v_0$ from the previous subsection, but now also assume given $\Xi_0\in E^{\mathbb{Z}\cap NI^\circ}$, the initial configuration of individual-channel states in the $N$th stochastic model.

DEFINITION 2.2. Suppose that $(\Omega,\mathcal{F},(\mathcal{F}_t)_{0\le t\le T},\mathbb{P})$ is a filtered probability space satisfying the usual conditions. Given a pair $(\mathbf{V}_t,\mathbf{\Xi}_t)_{0\le t\le T}$ of càdlàg adapted stochastic processes such that each sample path of $\mathbf{V}$ is a continuous map $[0,T]\to H_0^1(I)$ and $\mathbf{\Xi}_t$ is in $E^{\mathbb{Z}\cap NI^\circ}$ for all $t\in[0,T]$, we will say that they *satisfy the $N$th stochastic Hodgkin–Huxley equations* (S$_N$) *with initial conditions* $v_0,\Xi_0$ if

- (Regularity) The map $t\mapsto\frac{d}{dt}\mathbf{V}_t$ lies in $L^2_{H^{-1}(I)}[0,T]$ almost surely;
- (Dynamics: PDE)

$$\frac{d}{dt}\mathbf{V}_t = \Delta\mathbf{V}_t + \frac{1}{N}\sum_{i\in\mathbb{Z}\cap NI^\circ} c_{\mathbf{\Xi}_t(i)}(v_{\mathbf{\Xi}_t(i)} - \mathbf{V}_t(i/N))\delta_{i/N}$$

$$\forall t\in[0,T],\ \mathbb{P}\text{-a.s.}$$

  [we will refer to this equation as (S$_N$-PDE)];

- (Dynamics: jump)

$$\mathbb{P}(\boldsymbol{\Xi}_{t+h}(i) = \zeta \mid \boldsymbol{\Xi}_t(i) = \xi) = \alpha_{\xi,\zeta}(\mathbf{V}_t(i/N))h + \mathrm{o}_{h\downarrow 0}(h)$$
$$\forall t \in [0, T), h \in (0, T - t],$$

  with the coordinate processes $(\boldsymbol{\Xi}_{t+h}(i))_{h>0}$ independent to first order in $h$ as $h \downarrow 0$ conditional on $\mathcal{F}_t$ [we will refer to this system of equations and conditions as $(S_N$-jump)];

- (Initial conditions: PDE)

$$\mathbf{V}_0 = v_0;$$

- (Initial conditions: jump)

$$\boldsymbol{\Xi}_0 = \Xi_0;$$

- (Boundary conditions: PDE only)

$$\mathbf{V}_t(\pm \ell) = 0 \qquad \forall t \in [0, T].$$

2.4. *The goal of this paper.* Before we can state the main result of this paper we need a little more notation. For $\xi \in E$ we write $\mathbf{C}_{\xi,N}$ for the map $E^{\mathbb{Z} \cap N I^\circ} \to H^{-1}(I)$ given by

$$\mathbf{C}_{\xi,N}(\Xi) = \frac{1}{N} \sum_{i \in \mathbb{Z} \cap N I^\circ, \; \Xi(i) = \xi} \delta_{i/N}$$

[the Dirac deltas $\delta_{i/N}$ are readily interpreted as elements of $H^{-1}(I)$]; so $\mathbf{C}_{\xi,N}(\Xi)$ places a mass of $1/N$ on each point $i/N \in I^\circ$ at which $\Xi$ is in state $\xi$. We refer to it as the *empirical distribution for $\xi$*. We introduce the distributions $\mathbf{C}_{\xi,N}$ for each individual state $\xi$ to meet the notational needs of the subsequent analysis.

We are now ready to state the result:

THEOREM 2.3. *Let $\varepsilon > 0$, and suppose given initial conditions $v_0$, $p_{\xi,0}$. Then for any $N$ sufficiently large, say $N \geq N_1$, there exists an initial condition $\Xi_0$ for $(S_N)$ so that there is some "high-probability" $\Omega_1 \subseteq \Omega$ with $\mathbb{P}(\Omega \setminus \Omega_1) < \varepsilon$ and such that*

$$\sup_{0 \leq t \leq T} \|\mathbf{V}_t^{(N)} - \mathbf{v}_t\|_{H_0^1(I)} < \varepsilon,$$

$$\sup_{0 \leq t \leq T} \|\mathbf{C}_{\xi,N}(\boldsymbol{\Xi}_t^{(N)}) - \mathbf{p}_{\xi,t}\|_{H^{-1}(I)} < \varepsilon,$$

*on $\Omega_1$.*

Colloquially, this theorem tells us that as $N \to \infty$ the stochastic ion-channel model of the axon gives a time-evolution of the potential difference along the axon that converges to that given by the deterministic model, uniformly up to a given finite time horizon, in probability.

This theorem will be proved in Section 4. The overarching idea when proving theorems of this sort is often to find an inequality that bounds the growth rate of the deviation in terms of the values the deviation has taken thus far; or, by integrating this inequality, to bound the current value of the deviation in terms of some average of the values it has taken so far. The most common formalization of this idea, and the one we will rely on, is Gronwall's lemma; we quickly recall this here:

PROPOSITION 2.4 (Gronwall's lemma). *Suppose $T > 0$ and $f : [0, T] \to \mathbb{R}$ is continuous. Suppose further that there are constants $A, B > 0$ such that*

$$f(t) \le A + B \int_0^t f(s) \, ds$$

*for all $t \in [0, T]$. Then $f(t) \le Ae^{Bt}$ for all $t \in [0, T]$.*

2.5. *The three scales of the models.* One interesting feature of the stochastic Hodgkin–Huxley model is that it relates behavior on three distinct scales: the flow of charge at the scale of individual ions; the opening and shutting of ion channels at the scale of large protein molecules; and the working of the whole axon.

The stochastic model of the ion channel is faithful at the second and third of these, but uses a simplified description of the behavior at the first—the smallest—as a continuum charge distribution. As is standard, the "random" movement of a rarefied distribution of very many very small particles in a suitable medium (here ions in solution) is simplified to a continuum evolving in time according to a parabolic PDE.

In this sense the stochastic model "averages away" the random behavior at the smallest scale. However, it retains a detailed description of the intermediate scale, whereas the deterministic equation takes an average here also: the simultaneous states of a great many small ion channels are forgotten, with only the proportion of channels in each state within each small length of the axon being retained.

Thus we can think of the difference between the stochastic and deterministic models as one of resolution: although neither model can "see" the individual ions, the stochastic model can see single channels, whereas even these are beyond the deterministic model. In this sense the main result of this paper is a check that if we average out over the smallest scale to obtain the stochastic model, and then consider a suitable limit of this to represent the vanishing size of the intermediate scale, we recover the model obtained by averaging over both smallest and intermediate scales from the start.

Some slightly unusual features of the stochastic model can be traced back to this three-scale property of the system under study. More common applications of Markov processes to the modeling of real-world systems need consider only two scales, often corresponding to the smallest and largest of the above. In this simpler case the state of the full Markov process will typically describe the complete state of the system in terms of a (large) discrete collection of components, possibly distributed in space; this may then have a continuum limit (often deterministic, but sometimes still stochastic, depending on the regime) in which the small scale has undergone averaging and so only a single, large-scale picture remains. The many stochastic components in the full model are traded in for a more complicated large-scale description, often based on spatially variable quantities evolving following a PDE.

However, because we obtain our stochastic model by performing only some of the possible averaging, not all, we are left with *both* a large number of stochastic components (the ion channels), *and* a complicated, spatially variable deterministic system following a PDE (the membrane potential along the axon) coupled to them. We will find that this occasionally puts the analysis of the stochastic model slightly beyond the reach of more routine results in either stochastic processes or PDE, and so just a little thought is needed to combine both disciplines and obtain useful results. We will see this first when proving existence for the stochastic processes, and again when we come to the estimates for the growth of various related stochastic processes that we need for proving convergence.

**3. Preliminary results.**  In the first two subsections below we discuss various general facts about a relevant diffusion semigroup and about the sample paths of finite-state-space Markov processes. We then move on to discuss existence and regularity for our equations (D) and $(S_N)$.

3.1. *A diffusion semigroup.*  To prove our main existence and regularity results later in this section we will need some facts about the Feller semigroup $(P_t)_{t \geq 0}$ corresponding to Brownian motion in $I$ absorbed at the end-points of $I$.

This semigroup can help us because of the connection between diffusive PDE and Feller diffusion processes that allows us to rewrite (D-PDE) in the integral form

$$\mathbf{v}_t = P_t v_0 + \int_0^t P_{t-s} \left( \sum_{\xi \in E} c_\xi \mathbf{p}_{\xi,s} \cdot (v_\xi - \mathbf{v}_s) \right) ds,$$

and similarly for (S$_N$-PDE). In order to use this integral representation we first need to prove certain regularity properties of the semigroup.

LEMMA 3.1. *Let $y$ be in the interior of $I$. Then $P_t\delta_y$ is a smooth function on $I$ vanishing at the end-points for any $t > 0$. Furthermore:*

1. *there is some constant $C_1 > 0$, depending on $T$ but otherwise not on $t \in [0, T]$, such that for any continuous function $f : [0, t] \to \mathbb{R}$ we have that the function*

$$I \to \mathbb{R} : x \mapsto \int_0^t f(s) P_{t-s}\delta_y(x)\, ds$$

*is in $H_0^1(I)$ and satisfies the estimate*

$$\left\| \int_0^t f(s) P_{t-s}\delta_y(\cdot)\, ds \right\|_{H_0^1(I)} \leq C_1 \|f\|_\infty;$$

2. *for any fixed $\varepsilon > 0$ there is some constant $C_2(\varepsilon)$, depending on $\varepsilon$ and $T$ but otherwise not on $t \in [0, T]$, such that for any continuous function $f : [0, t] \to \mathbb{R}$ we have*

$$\left\| \int_0^{t-\varepsilon} f(s) P_{t-s}\delta_y(\cdot)\, ds \right\|_{H_0^1(I)} \leq C_2(\varepsilon) \int_0^{t-\varepsilon} |f(s)|\, ds.$$

*for any $t \in [0, T]$.*

REMARK. Note that the imposition of a fixed $\varepsilon > 0$ in the second estimate is necessary; without it the result can be made to fail for any given choice of $C_2$ by choosing $f$ to be zero apart from in $(t - \eta, t)$, where it rapidly becomes very large, for some sufficiently small $\eta > 0$.

PROOF OF LEMMA 3.1. By additivity we may assume $f \geq 0$. One-dimensional Brownian motion has the transition density

$$p_t(x, y) = \frac{1}{\sqrt{2\pi t}} e^{-|x-y|^2/2t},$$

so that for $f \in C_{\mathrm{b}}(\mathbb{R})$ the corresponding Feller semigroup $(Q_t)_{t \geq 0}$ is given by

$$Q_t f(x) = \int_{\mathbb{R}} f(y) \frac{1}{\sqrt{2\pi t}} e^{-|x-y|^2/2t}\, \mu(dy).$$

Now our semigroup $(P_t)_{t \geq 0}$ corresponds to Brownian motion absorbed at the end-points of $I$ (see, e.g., Chapter 24 of [14]). This semigroup has the modified transition density

$$p_t^I(x, y) = p_t(x, y) - \mathbb{E}_x(p_{t-\tau}(W_\tau, y)\chi_{\{\tau < t\}})$$
$$= p_t(x, y) - \mathbb{E}_x(p_{t-\tau}(\ell, y)\chi_{\{\tau < t, W_\tau = \ell\}})$$
$$- \mathbb{E}_x(p_{t-\tau}(-\ell, y)\chi_{\{\tau < t, W_\tau = -\ell\}}),$$

where $W$ is our Brownian motion and $\tau$ is the hitting time of the boundary of $I$. Applying $P_t$ to the Dirac point-mass $\delta_y$ we recover precisely this expression for $p_t^I(x, y)$. That $P_t \delta_y$ is a smooth function vanishing at the end-points of $I$ now follows at once.

It remains to establish the two estimates. From the above we have

$$\int_0^t f(s) P_{t-s} \delta_y(x)\, ds = \int_0^t f(s) \frac{1}{\sqrt{2\pi(t-s)}} e^{-|x-y|^2/2(t-s)}\, ds$$
$$- \int_0^t f(s) \mathbb{E}_x(p_{t-\tau}(\ell, y)\chi_{\{\tau < t-s, W_\tau = \ell\}})\, ds$$
$$- \int_0^t f(s) \mathbb{E}_x(p_{t-s-\tau}(-\ell, y)\chi_{\{\tau < t-s, W_\tau = -\ell\}})\, ds.$$

It suffices to prove the desired regularity for each of these three integrals separately. That both of the estimates hold with suitable constants (the second depending on $\varepsilon$) for the second and third integrals is clear, since $y$ is fixed away from $\pm\ell$ and so the expressions inside the expectations $\mathbb{E}_x$ are uniformly bounded functions of $x$ with uniformly bounded equicontinuous derivatives as $s$ varies in $(0, t)$.

We are left with the first integral, which we break into two pieces.

First we estimate the integral over $(t-\varepsilon, t)$. We have

$$\left| \int_{t-\varepsilon}^t f(s) \frac{1}{\sqrt{2\pi(t-s)}} e^{-|x-y|^2/2(t-s)}\, ds \right|$$
$$\leq \|f\|_\infty \int_{t-\varepsilon}^t \frac{1}{\sqrt{2\pi(t-s)}} e^{-|x-y|^2/2(t-s)}\, ds,$$

and now, making the substitution $s = t - 1/u^2$, this is

$$\int_{1/\sqrt{\varepsilon}}^\infty \frac{1}{\sqrt{2\pi} u^2} e^{-(u(x-y))^2/2}\, du.$$

For any given $\varepsilon$ this is clearly a smooth function of $x$, since the integrand over $(1/\sqrt{\varepsilon}, \infty)$ is dominated by $1/u^2$. It is also clear that it converges to 0 in $\|\cdot\|_{L^2(I)}$, by dominated convergence; by nonnegativity the same is true of the original integral involving $f$.

We may also differentiate our integral over $(t-\varepsilon, t)$ with respect to $x$ under the integral sign (except, possibly, at $x = y$) to obtain the new integral

$$-\int_{t-\varepsilon}^{t} f(s) \frac{x-y}{\sqrt{2\pi}(t-s)^{3/2}} e^{-|x-y|^2/2(t-s)}\, ds,$$

which is bounded in absolute value by

$$\|f\|_{\infty} \int_{t-\varepsilon}^{t} \frac{|x-y|}{\sqrt{2\pi}(t-s)^{3/2}} e^{-|x-y|^2/2(t-s)}\, ds.$$

Using the substitution $s = t - 1/u^2$ again this becomes

$$\int_{1/\sqrt{\varepsilon}}^{\infty} \frac{|x-y|}{\sqrt{2\pi}} e^{-(u(x-y))^2/2}\, du.$$

Making the second substitution $u = w/|x-y|$, this becomes in turn

$$\|f\|_{\infty} \frac{1}{\sqrt{2\pi}} \int_{|x-y|/\sqrt{\varepsilon}}^{\infty} e^{-w^2/2}\, dw$$

[noting the cancelation of two factors of $(x-y)$], which is bounded as $x \to y$ and so is also in $L^2(I)$ as a function of $x$. Another appeal to the Dominated Convergence Theorem completes the proof that this tends to zero in $\|\cdot\|_{L^2(I)}$ as $\varepsilon \to 0$.

Therefore for any $\eta > 0$ we can choose $\varepsilon > 0$ so small that

$$\left\| \int_{t-\varepsilon}^{t} f(s) P_{t-s} \delta_y(\cdot)\, ds \right\|_{H_0^1(I)}$$

$$\leq \left\| \int_{t-\varepsilon}^{t} f(s) P_{t-s} \delta_y(\cdot)\, ds \right\|_{L^2(I)} + \left\| D\left( \int_{t-\varepsilon}^{t} f(s) P_{t-s} \delta_y(\cdot)\, ds \right) \right\|_{L^2(I)}$$

$$\leq \|f\|_{\infty} \left( \left\| \int_{1/\sqrt{\varepsilon}}^{\infty} \frac{1}{\sqrt{2\pi} u^2} e^{-u^2 |(\cdot)-y|^2/2}\, du \right\|_{L^2(I)}$$

$$+ \left\| \frac{1}{\sqrt{2\pi}} \int_{|(\cdot)-y|/\sqrt{\varepsilon}}^{\infty} e^{-w^2/2}\, dw \right\|_{L^2(I)} \right)$$

$$< \frac{1}{2} \eta \|f\|_{\infty}.$$

Having done so, the uniform boundedness and smoothness properties of $P_{t-s}\delta_y(x)$ (as a function of $x$) for $s$ bounded away from $t$ give at once some $\delta > 0$ such that if $t \leq T$ and $\|f\|_{L^1[0,T]} \leq \delta$, then

$$\left\| \int_{0}^{t-\varepsilon} f(s) P_{t-s} \delta_y(\cdot)\, ds \right\|_{H_0^1(I)} < \frac{1}{2} \eta.$$

For suitable $\eta$ this is just the second estimate that we wanted, and adding to our previous inequality gives also the first estimate. $\qquad\square$

3.2. *Finite-state-space Markov processes.* We recall here some facts about Markov processes with a finite state space $S$ and the corresponding space of sample paths $D_S[0, \infty)$; or rather, as we will need, the space $D_S[0, T]$ of paths taken only up to (and at) a finite time-horizon $T$. This space of paths is treated comprehensively in Chapter 3 of Ethier and Kurtz [5].

For any complete, separable metric space $S$ the space $D_S[0, \infty)$ of càdlàg paths from $[0, \infty)$ into $S$ is also complete and separable when endowed with its Skorohod topology; in particular, this is so if $S$ is a finite set with its discrete metric. The same is clearly true of our space $D_S[0, T]$. We will write $\pi_t$ for the time-$t$ projection map $D_S[0, \infty) \to S : \omega \mapsto \omega_t$. The process $(\pi_t)_{t \geq 0}$ is referred to as the *canonical process*, and defines the *canonical filtration*

$$\mathcal{F}_t = \sigma(\{\pi_s : s \leq t\}).$$

In the case of $S$ a finite set each path $\omega \in D_S[0, T]$ is completely characterized by the following data:

- the total number $N(\omega)$ of jumps performed by the path (this is always finite);
- the sequence of numbers $\sigma_j(\omega) > 0$, $j = 1, 2, \ldots, N(\omega)$, giving the time of the $j$th jump (for convenience we also set $\sigma_0 = 0, \sigma_{N(\omega)+1} = T$);
- the sequence of states $\xi_j$, $j = 0, 1, \ldots, N(\omega)$, giving the starting state for $j = 0$ and the landing state $\omega_{\sigma_j(\omega)}$ after the $j$th jump for $j \geq 1$.

It is clear that all of the above define measurable (in fact, continuous) functions on $D_S[0, T]$.

When we come to construct a solution to our stochastic equations, we will need the following explicit computation of the absolutely continuous change of measure on the path space $D_S[0, T]$ implied by the Girsanov theorem in the context of Markov processes with finite state spaces; in this sense it is an analogue of the Cameron–Martin theorem. The required theory and calculations are treated in Chapter III, Section 5 of Jacod and Shiryaev [13] (in a rather more general setting).

Lemma 3.2. *Suppose that $\mathbb{P}_1$ is a probability measure on $\Omega = D_S[0, T]$ for which the canonical process $(\pi_t)_{t \in [0,T]}$ is a Markov process with all jump rates equal to 1; and suppose also that for each $\xi, \zeta \in S$, $\xi \neq \zeta$, we are given a progressively measurable function $\lambda_{\xi,\zeta} : \Omega \times [0, T] \to [0, \infty)$ with $\lambda_{\xi,\zeta}(\omega, \cdot)$ continuous for every $\omega$. Define $h : \Omega \to \mathbb{R}$ by*

$$h(\omega) = \frac{\prod_{j=0}^{N(\omega)} \exp(-\int_{\sigma_j(\omega)}^{\sigma_{j+1}(\omega)} \lambda_{\xi_j(\omega)}(\omega, s)\, ds) \lambda_{\xi_j(\omega), \xi_{j+1}(\omega)}(\omega, \sigma_{j+1}(\omega))}{\prod_{j=0}^{N(\omega)} \exp(-(\sigma_{j+1}(\omega) - \sigma_j(\omega)))}$$

$$= e^T \prod_{j=0}^{N(\omega)} \exp\left(-\int_{\sigma_j(\omega)}^{\sigma_{j+1}(\omega)} \lambda_{\xi_j(\omega)}(\omega, s)\, ds\right) \lambda_{\xi_j(\omega), \xi_{j+1}(\omega)}(\omega, \sigma_{j+1}(\omega))$$

*and* $\mathbb{P} = \mathbb{P}_1 \llcorner h$. *Then under* $\mathbb{P}$ *the canonical process is a time-inhomogeneous Markov process with jump rates* $\lambda_{\xi,\zeta}(\omega, t)$ *for* $t \in [0, T]$, *and* $\mathbb{P}$ *is the unique probability on* $\Omega$ *with this property.*

### 3.3. *Existence and uniqueness for the deterministic equations.*

PROPOSITION 3.3. *There is a unique weak solution to* (D), *and it satisfies* $v_- \leq \mathbf{v}_t \leq v_+$ *for all* $t$.

PROOF. This is a classical example of the use of fixed-point theorems and Gronwall's lemma in the study of nonlinear parabolic PDE, and we will not describe it here (we will see a more complicated example in the existence and uniqueness result for the stochastic equations anyway). A thorough and readable treatment is Lamberti's [16], although the first existence and uniqueness results are for *strong* solutions and go back to Evans and Shenk [6]. Note that both of these papers give an analysis specific to the original Hodgkin–Huxley equations, with particular forms for the states $\xi$ and proportions $\mathbf{p}_\xi$; the method of analysis, however, extends to our case immediately. □

### 3.4. *Existence and uniqueness for the stochastic equations.* Concerning the stochastic equations, our later proof of convergence will need only a suitable form of weak existence of solutions, and not uniqueness. However, it seems only natural to include proofs of both existence and uniqueness here.

In constructing the process $(\mathbf{V}, \mathbf{\Xi})$ we will need to introduce a particular underlying filtered probability space $(\Omega, \mathcal{F}, (\mathcal{F}_t)_{0 \leq t \leq T}, \mathbb{P})$, even though our later convergence results hold for a suitable process on any such space. We will choose the probability space with some additional structure that allows us to interpret any $\omega \in \Omega$ as a driving signal from which we can (almost surely) construct a corresponding sample path of our desired process. Thus, the choice of a particular filtered probability space, and the subsequent construction of a probability on it, can be thought of as the choice of how to mimic the randomness apparent in the real-world system.

Now, the process we want takes values in the overall state space $H_0^1(I) \times E^{\mathbb{Z} \cap NI^\circ}$, which is far from locally compact, and so the above construction is not contained within the standard machinery of Feller processes: we will need to construct our process with a little more care. We remark that this lack of local compactness is an artifact of our model's two different "small" scales (see Section 2.5): the function space arises as a result of the averaging over the "very small," and so leaves us to cope with the infinite-dimensional topology of that function space, while we still want to model the "fairly small" scale stochastically.

The key to our construction is to observe from the dynamics of $(S_N)$ that, if we already knew the final form of the process $\boldsymbol{\Xi}$, then for a fixed sample path of $\boldsymbol{\Xi}$ the evolution of the corresponding sample path of the process $\mathbf{V}$ would be deterministic. Happily, regarded just as a PDE, $(S_N\text{-PDE})$ is a nonlinear parabolic equation, and at the core of the theory of these is a standard procedure for proving existence of solutions. It relies on a handful of classical fixed-point theorems for Banach spaces and a few ways of choosing how to apply one of them to a suitable function space. We will use a slight modification of that procedure applied $\omega$-by-$\omega$.

To do this, we again convert the PDE to an integral equation:

$$\mathbf{V}_t = P_t v_0 + \int_0^t \frac{1}{N} P_{t-s} \left( \sum_{i \in \mathbb{Z} \cap N I^\circ} c_{\boldsymbol{\Xi}_s} (v_{\boldsymbol{\Xi}_s} - \mathbf{V}_s) \delta_{i/N} \right) ds$$

$$= P_t v_0 + \frac{1}{N} \sum_{i \in \mathbb{Z} \cap N I^\circ} \int_0^t c_{\boldsymbol{\Xi}_s} (v_{\boldsymbol{\Xi}_s} - \mathbf{V}_s)(P_{t-s} \delta_{i/N}) \, ds,$$

where $(P_t)_{t \geq 0}$ is the Feller semigroup from Section 3.1.

We will apply our chosen fixed-point result to the Banach space $X = C_{H_0^1(I)}[0,T]$; the required differentiability properties of the function $\mathbf{V}$ will then follow from the integral equation. It is worth commenting on this choice of space. For many applications in PDE the larger space $C_{L^2(I)}[0,T]$, with its less restrictive topology, would be the appropriate choice. We are forced to work with the smaller space by the slightly unusual nature of our PDE: the right-hand term contains a linear combination of Dirac measures, and so is not itself a function but only a member of $H^{-1}(I)$. It will turn out when we construct our map from $X$ to itself below that the smoothing properties of the heat semigroup are not enough to give a continuous self-map of $C_{L^2(I)}[0,T]$, and so we are forced to work with the more complicated space.

Motivated by these observations, we will attempt to construct $(\mathbf{V}, \boldsymbol{\Xi})$ on the space $\Omega = D_{E^{\mathbb{Z} \cap N I^\circ}}[0,T]$ of càdlàg paths from $[0,T]$ into $E^{\mathbb{Z} \cap N I^\circ}$, regarded as a filtered space as described in Section 3.2. We will start by defining a "simple" probability $\mathbb{P}_1$ on this space, and will later obtain the desired probability $\mathbb{P}$ as a suitable indefinite-integral measure with respect to $\mathbb{P}_1$, using Lemma 3.2. This use of the path space makes more concrete the above-mentioned driving signal interpretation of $\omega \in \Omega$.

We choose our probability $\mathbb{P}_1$ by specifying that under it the canonical process $(\pi_t)_{t \in [0,T]}$ is a Feller process with all jump rates between different configurations in $E^{\mathbb{Z} \cap N I^\circ}$ equal to 1 (the existence of such a probability is guaranteed by the usual theory of Feller processes).

PROPOSITION 3.4. *Fix $N \geq 1$ and let the filtered space $(\Omega, \mathcal{F}, (\mathcal{F}_t)_{0 \leq t \leq T}, \mathbb{P})$ be as above. Then the space carries a solution to $(S_N)$ with initial conditions $v_0$ and $\Xi_0$, and the law of this solution is unique.*

PROOF. The proof strategy is as follows:

1. Obtain from the integral form of $(S_N\text{-PDE})$ an equation for the trajectory $(\mathbf{V}_t(\omega))_{0 \leq t \leq T}$ for a given input signal $\omega \in D_{E^{\mathbb{Z} \cap NI^\circ}}[0, T]$. This equation will then be of the form $\mathbf{V}(\omega) = \Psi(\mathbf{V}(\omega), \omega)$ for a suitable map $\Psi \colon X \times \Omega \to X$.

2. Show that for each $\omega$ separately the map $\Psi_\omega = \Psi(\cdot, \omega) \colon X \to X$:
   (a) is continuous;
   (b) has the compactness property needed for Schaefer's Fixed-Point Theorem (the appropriate fixed-point theorem for this proof, since we do not have any obvious contraction mapping);
   (c) has the boundedness property needed for Schaefer's theorem;
   and so obtain a nonempty set of fixed points for every $\omega \in \Omega$. It is here that we will need our preliminary lemmas about the semigroup $(P_t)_{t \geq 0}$. At this point we will also show that the trajectory $\mathbf{V}(\omega)$ is unique, given $\omega$.

3. Show that $\Psi \colon X \times \Omega \to X$ is measurable, and hence that the set of pointwise fixed points $\{(\mathbf{U}, \omega) \colon \Psi(\mathbf{U}, \omega) = \mathbf{U}\}$ is measurable and has a nonempty section above every $\omega \in \Omega$, and apply the Measurable Selector Theorem to give the function $\mathbf{V}$.

4. Having thus obtained a suitable trajectory $\mathbf{V}(\omega)$ for each $\omega$, varying progressively measurably with the sample path $\omega$ in $D_{E^{\mathbb{Z} \cap NI^\circ}}[0, T]$, we can use Lemma 3.2 to write down a suitable Radon–Nikodým derivative for a new probability $\mathbb{P}$ with respect to $\mathbb{P}_1$ so that under $\mathbb{P}$ the equation $(S_N\text{-jump})$ is also satisfied.

5. Finally, uniqueness will follow from the uniqueness for the PDE corresponding to a single $\omega$ proved at the end of Step 2, and the uniqueness part of Lemma 3.2.

STEP 1. For each $\omega \in D_{E^{\mathbb{Z} \cap NI^\circ}}[0, T]$ we need to solve the integral equation

$$\mathbf{V}_t(\omega) = P_t v_0 + \frac{1}{N} \sum_{i \in \mathbb{Z} \cap NI^\circ} \int_0^t c_{\omega_s(i)}(v_{\omega_s(i)} - \mathbf{V}_s(\omega)(i/N))(P_{t-s}\delta_{i/N})\, ds;$$

we let $\Psi(\mathbf{V}(\omega), \omega)_t$ be the expression on the right-hand side.

STEP 2. Fix $\omega$. Recall Schaefer's Fixed-Point Theorem (see, e.g., Theorem 4 of Section 9.2.2 in Evans [7]):

SCHAEFER'S FIXED-POINT THEOREM. *Suppose that $X$ is a Banach space and that $\Psi : X \to X$ is a continuous map that converts bounded sequences to precompact sequences. Assume further that the set*

$$\{u \in X : u = \lambda \Psi u \text{ for some } \lambda \in [0,1]\}$$

*is bounded in $X$. Then $\Psi$ has a fixed point.*

We will check these relevant properties separately for $\Psi = \Psi_\omega$. We write $\Psi_\omega$ as

$$\Psi_\omega(\mathbf{U}) = P_t v_0 + \frac{1}{N} \sum_{i \in \mathbb{Z} \cap NI^\circ} \Psi_{\omega,i}(\mathbf{U})$$

where

$$\Psi_{\omega,i}(\mathbf{U})_t = \int_0^t c_{\omega_s(i)} (v_{\omega_s(i)} - \mathbf{U}_s(i/N))(P_{t-s}\delta_{i/N}) \, ds.$$

STEP 1 (Continuity). It suffices to prove this for each $\Psi_{\omega,i}$, for which it will follow from Lemma 3.1. Suppose $\mathbf{U}, \mathbf{W} \in X$. Then

$$\Psi_{\omega,i}(\mathbf{U})_t - \Psi_{\omega,i}(\mathbf{W})_t$$
$$= \int_0^t (c_{\omega_s(i)}(v_{\omega_s(i)} - \mathbf{U}_s(i/N))$$
$$\quad - c_{\omega_s(i)}(v_{\omega_s(i)} - \mathbf{W}_s(i/N)))(P_{t-s}\delta_{i/N}) \, ds$$
$$= \int_0^t c_{\omega_s(i)}(\mathbf{W}_s(i/N) - \mathbf{U}_s(i/N))(P_{t-s}\delta_{i/N}) \, ds.$$

Since by Poincaré's inequality the norm

$$\sup_{0 \le t \le T} \|\mathbf{U}_t - \mathbf{W}_t\|_{H_0^1(I)}$$

is stronger than

$$\sup_{0 \le t \le T} \|\mathbf{U}_t - \mathbf{W}_t\|_\infty$$

to within a multiplicative constant, it follows that by selecting the former sufficiently small we may make the multiplier $c_{\omega_s(i)}(v_{\omega_s(i)} - \mathbf{W}_s(i/N))$ of $P_{t-s}\delta_{i/N}$ uniformly small in the above integrand, and so make the norm

$$\|\Psi_{\omega,i}(\mathbf{U})_t - \Psi_{\omega,i}(\mathbf{W})_t\|_{H_0^1(I)}$$

as small as we please uniformly in $t \in [0,T]$, by Lemma 3.1.

STEP 2 (Compactness). Here also it suffices to consider each $\Psi_{\omega,i}$ separately. Compactness now follows directly from the estimates of Lemma 3.1, which allow us to approximate the integral expression by a linear combination of the functions $P_s \delta_{i/N}$ for $s$ taken from some sufficiently large finite subset of $(0, t)$.

STEP 3 (Boundedness). Suppose that $\mathbf{U} \in X$ has $\mathbf{U} = \lambda \Psi_\omega(\mathbf{U})$ for some $\lambda \in [0, 1]$. Writing this out more fully, it reads

$$\mathbf{U}_t = \lambda P_t v_0 + \lambda \frac{1}{N} \sum_{i \in \mathbb{Z} \cap NI^\circ} \int_0^t c_{\omega_s(i)}(v_{\omega_s(i)} - \mathbf{U}_s(i/N))(P_{t-s}\delta_{i/N}) \, ds.$$

Hence

$$
\begin{aligned}
\|\mathbf{U}_t\|_{H_0^1(I)} & \\
\leq \lambda & \|P_t v_0\|_{H_0^1(I)} \\
& + \lambda \frac{1}{N} \sum_{i \in \mathbb{Z} \cap NI^\circ} \left\| \int_0^t c_{\omega_s(i)}(v_{\omega_s(i)} - \mathbf{U}_s(i/N))(P_{t-s}\delta_{i/N}) \, ds \right\|_{H_0^1(I)} \\
\leq \lambda & \|P_t v_0\|_{H_0^1(I)} + \lambda \frac{1}{N} \sum_{i \in \mathbb{Z} \cap NI^\circ} \left\| \int_0^t c_{\omega_s(i)} v_{\omega_s(i)}(P_{t-s}\delta_{i/N}) \, ds \right\|_{H_0^1(I)} \\
& + \lambda \frac{1}{N} \sum_{i \in \mathbb{Z} \cap NI^\circ} \left\| \int_0^t c_{\omega_s(i)} \mathbf{U}_s(i/N)(P_{t-s}\delta_{i/N}) \, ds \right\|_{H_0^1(I)}.
\end{aligned}
$$

Rather than try to bound the growth of $\|\mathbf{U}_t\|_{H_0^1(I)}$ directly using the above inequality, we consider the maximal function $u(t) = \max_{0 \leq s \leq t} \|\mathbf{U}_t\|_{H_0^1(I)}$. By Lemma 3.1, the first of the above sums is bounded by a fixed constant (since there are only finitely many possible values for $c_{\xi,\zeta}$ and $v_\xi$), which can be chosen independent of $t \in [0, T]$. For the second sum, we break the integral in each term into two pieces: an integral over a small interval $(t - \varepsilon, t)$ that we can bound using the length $\varepsilon$ of the interval, and another over what remains that we can bound because $P_{t-s}\delta_{i/N}$ is more regular there. This idea is similar to that in the proof of Lemma 3.1. We select $\varepsilon$ so small that

$$\max_{\xi \in E} |c_\xi| \left\| \int_{t-\varepsilon}^t P_{t-s}\delta_{i/N} \, ds \right\|_{H_0^1(I)} \leq \frac{1}{2C_1},$$

and so, using Lemma 3.1,

$$\left\| \int_{t-\varepsilon}^t c_{\omega_s(i)} \mathbf{U}_s(i/N)(P_{t-s}\delta_{i/N}) \, ds \right\|_{H_0^1(I)} \leq \tfrac{1}{2} u(t).$$

Now the second estimate of Lemma 3.1 gives

$$\left\| \int_0^{t-\varepsilon} c_{\omega_s(i)} \mathbf{U}_s(i/N)(P_{t-s}\delta_{i/N}) \, ds \right\|_{H_0^1(I)} \le C_2(\varepsilon) \max_{\xi \in E} |c_\xi| \int_0^{t-\varepsilon} \|\mathbf{U}_s\|_\infty \, ds$$

$$\le C C_2(\varepsilon) \max_{\xi \in E} |c_\xi| \int_0^t u(s) \, ds,$$

where $C$ is the constant so that $\|\cdot\|_\infty \le C \|\cdot\|_{H_0^1(I)}$, guaranteed by Poincaré's inequality. Reassembling the above inequalities, we find

$$\|\mathbf{U}_t\|_{H_0^1(I)} \le A + \lambda B \int_0^t u(s) \, ds + \frac{\lambda}{2} u(t),$$

where

$$A = \lambda \|P_t v_0\|_{H_0^1(I)} + \lambda \frac{1}{N} \sum_{i \in \mathbb{Z} \cap NI^\circ} \sup_{0 \le t \le T} \left\| \int_0^t c_{\omega_s(i)} v_{\omega_s(i)} (P_{t-s}\delta_{i/N}) \, ds \right\|_{H_0^1(I)}$$

and

$$B = C C_2(\varepsilon) \max_{\xi \in E} |c_\xi|,$$

and so

$$u(t) \le 2A + 2B \int_0^t u(s) \, ds$$

(as $\lambda \le 1$). Now Gronwall's lemma gives a fixed bound on $u(t)$ for $t \in [0, T]$, and so we have the desired bound in the space $X$.

COMPLETION OF STEP 2. Thus the set $Y_\omega$ of fixed points of $\Psi_\omega$ is nonempty for every $\omega \in \Omega$.

We also wish to show that the solution to our PDE (or, equivalently, integral equation) for a fixed $\omega$ is unique, that is, that $Y_\omega$ is a singleton for every $\omega \in \Omega$. This is an easy calculation; for suppose $\mathbf{V}^j(\omega) = \Psi_\omega(\mathbf{V}^j)$, $j = 1, 2$, are two solutions for some $\omega \in \Omega$. Then, taking the difference, we find that $\mathbf{U} := \mathbf{V}^1(\omega) - \mathbf{V}^2(\omega)$ satisfies

$$\mathbf{U}_t = \frac{-1}{N} \sum_{i \in \mathbb{Z} \cap NI^\circ} \int_0^t c_{\omega_s(i)} \mathbf{U}_s(i/N)(P_{t-s}\delta_{i/N}) \, ds.$$

Repeating the trick of breaking $(0, t)$ into $(0, t - \varepsilon)$ and $(t - \varepsilon, t)$ and applying Lemma 3.1, as used in Step 3 above, we now obtain the estimate

$$\|\mathbf{U}_t\|_{H_0^1(I)} \le B \int_0^t u(s) \, ds + \frac{1}{2} u(t),$$

where once again $u(t) = \max_{0 \le s \le t} \|\mathbf{U}_t\|_{H_0^1(I)}$. Hence

$$u(t) \le B \int_0^t u(s)\, ds + \tfrac{1}{2} u(t),$$

for all $t \in [0, T]$, and by Gronwall's lemma $u \equiv 0$, and thus $\mathbf{V}^1(\omega) = \mathbf{V}^2(\omega)$.

STEP 3. Now set

$$Y = \{(\mathbf{U}, \omega) \in X \times \Omega \colon \mathbf{U} \text{ is a fixed point for } \Psi_\omega\}$$
$$= \{(\mathbf{U}, \omega) \in X \times \Omega \colon \Psi(\mathbf{U}, \omega) - \mathbf{U} = 0\}.$$

Since the map $X \times \Omega \to X$ sending $(\mathbf{U}, \omega)$ to $\Psi(\mathbf{U}, \omega) - \mathbf{U}$ is measurable, so is the set $Y$; and by Step 2, all of the sections $Y \cap (X \times \{\omega\}) = Y_\omega$ are nonempty. $X$ is a separable Banach space, and so the Measurable Selector Theorem (for a suitable version see, e.g., Theorem 10.1 in the Appendix to Ethier and Kurtz [5]) guarantees a measurable function $\mathbf{V} \colon \Omega \to X$ with $\mathbf{V}(\omega) \in Y_\omega$ for all $\omega$.

STEP 4. Now we can obtain $\mathbb{P}$ as a suitable indefinite-integral measure with respect to $\mathbb{P}_1$ so that $(S_N\text{-jump})$ is satisfied. This follows from Lemma 3.2 by taking the global jump-rates in $E^{\mathbb{Z} \cap NI^\circ}$ that correspond to the independent jump-rates of single components given by $(S_N\text{-jump})$; that is, for distinct configurations $\Xi_1, \Xi_2 \in E^{\mathbb{Z} \cap NI^\circ}$,

$$\lambda_{\Xi_1, \Xi_2}(\omega, t) = \begin{cases} 0, & \text{if } \Xi_1,\ \Xi_2 \text{ differ in two coordinates} \\ & \text{or more,} \\ \alpha_{\xi, \zeta}(\mathbf{V}_t(\omega)(i/N)), & \text{if } \Xi_1(i) = \xi \ne \zeta = \Xi_2(i) \\ & \text{and } \Xi_1,\ \Xi_2 \text{ agree everywhere else.} \end{cases}$$

In order to apply the lemma, we need only note that the sample paths of the process $\mathbf{V}_t$ are continuous and depend only on the evolution of the jump process up to time $t$ [from the form of the integral solution to $(S_N\text{-PDE})$], and so the rate processes $\lambda_{\Xi_1, \Xi_2}$ are continuous and progressively measurable.

STEP 5. It remains to prove uniqueness in law. Suppose now that $(\mathbf{V}_t^0, \Xi_t^0)_{0 \le t \le T}$ is some solution to $(S_N)$ on some abstract filtered space $(\Omega^0, \mathcal{F}^0, (\mathcal{F}_t^0)_{0 \le t \le T}, \mathbb{P}^0)$. Then, from $(S_N\text{-PDE})$, we know that

$$\mathbf{V}_t^0 = P_t v_0 + \frac{1}{N} \sum_{i \in \mathbb{Z} \cap NI^\circ} \int_0^t c_{\Xi_s^0(i)}(v_{\Xi_s^0(i)} - \mathbf{V}_s^0(i/N))(P_{t-s}\delta_{i/N})\, ds$$

for all $t \in [0, T]$, $\mathbb{P}^0$-almost surely. From the end of Step 2 above we know that this equation specifies $\mathbf{V}^0$ uniquely, given the sample path $\Xi^0$, and therefore $\mathbf{V}^0 = \mathbf{V}(\mathbf{\Xi}^0)$, $\mathbb{P}^0$-almost surely. It now follows from (S$_N$-jump) that the law of $\mathbf{\Xi}^0$ under $\mathbb{P}^0$ must have the same inhomogeneous Markov property as described for $\mathbb{P}$ in Step 4, and so, by the uniqueness part of Lemma 3.2, these probability measures must be equal. The result follows.  □

3.5.  *Three additional regularity results.*  In this subsection we collect three additional results that will be needed later.

The first is almost immediate.

PROPOSITION 3.5.  *Suppose* $v_0 \in H_0^1(I)$, $v_- \leq v_0 \leq v_+$, *and consider* (D) *with initial conditions* $(v_0, (p_{\xi,0})_{\xi \in E})$ *and each* (S$_N$) *with initial conditions* $(v_0, \Xi_0^{(N)})$. *Then there is some constant* $C_3 > 0$ *(independent of the exact initial conditions and of* $N$*) such that*

$$\|\mathbf{v}_t\|_\infty, \|\mathbf{V}_t^{(N)}\|_\infty \leq C_3$$

*for all* $N \geq 1$ *and all* $t \in [0, T]$.

PROOF.  It follows from Proposition 3.3 that the constant $\max_{\xi \in E} |v_\xi|$ itself works for the deterministic trajectory $\mathbf{v}$; it will therefore suffice to find a constant that works simultaneously for all of the stochastic equations. However, we have already proved this while establishing the boundedness in Step 3 in the proof of existence for the stochastic equations, for the constants $A$ and $B$ used there did not depend on $N$, $t \in [0, T]$ or $\omega$ and we may just set $\lambda = 1$. This completes the proof.  □

REMARK.  In fact it seems intuitively clear that $C_3 = \max_{\xi \in E} |v_\xi|$ should work for the stochastic equations also, but I have not proved this.

The next result tells us a little more about the regularity of the map $\mathbf{v}$:

LEMMA 3.6.  *Suppose* $v_0$ *and* $p_\xi$, $\xi \in E$, *have common Lipschitz constant* $K < \infty$, *and that* $\mathbf{v}$ *and* $\mathbf{p}_\xi$ *are a solution to* (D). *Then*

$$\sup_{0 \leq t \leq T} \|D\mathbf{v}_t\|_\infty < \infty.$$

REMARK.  An $L^\infty$ bound such as this seems a little odd for a PDE for which existence of solutions is most naturally studied in the Sobolev space $H_0^1(I)$, and which has a distinctly quadratic flavor owing to its diffusive nature. We will need this $L^\infty$ estimate later as a consequence of the very specific form of the nonlinear term in (D-PDE).

PROOF OF LEMMA 3.6. We observe first that since $\mathbf{v}\colon [0,T] \to H_0^1(I)$ is continuous, so is the function from $[0,T] \times I \to \mathbb{R}$ defined by $(t,y) \mapsto \mathbf{v}_t(y)$. Therefore we do have

$$\sup_{0 \le t \le T} \|\mathbf{v}_t\|_\infty < \infty,$$

and hence also

$$\sup_{0 \le t \le T} \left\| \sum_{\xi \in E} c_\xi \mathbf{p}_{\xi,t} \cdot (v_\xi - \mathbf{v}_t) \right\|_\infty < \infty.$$

Now the proof makes another use of conversion to an integral equation, this time for (D-PDE):

$$\mathbf{v}_t = P_t v_0 + \int_0^t P_{t-s} \left( \sum_{\xi \in E} c_\xi \mathbf{p}_{\xi,s} \cdot (v_\xi - \mathbf{v}_s) \right) ds$$

$$= P_t v_0 + \sum_{\xi \in E} c_\xi \int_0^t P_{t-s} (\mathbf{p}_{\xi,s} \cdot (v_\xi - \mathbf{v}_s)) \, ds.$$

The easiest way to proceed is to appeal to the specific form for $t > 0$ of the transition density $p_t^I(x,y)$ corresponding to the semigroup $(P_t)_{t \ge 0}$, as in the proof of Lemma 3.1:

$$p_t^I(x,y) = \frac{1}{\sqrt{2\pi t}} e^{-|x-y|^2/2t} - \mathbb{E}_x(p_{t-\tau}(\ell,y) \chi_{\{\tau < t, W_\tau = \ell\}})$$

$$- \mathbb{E}_x(p_{t-\tau}(-\ell,y) \chi_{\{\tau < t, W_\tau = -\ell\}}),$$

where $\tau$ is the hitting time of the boundary of $I$ for a standard Brownian motion.

Using this, we can write for a continuous function $f\colon I \to \mathbb{R}$ and $x \in \mathbb{R}$

$$P_t f(x) = \int_I f(y) \frac{1}{\sqrt{2\pi t}} e^{-|x-y|^2/2t} \, \mu(dy)$$

$$- \int_I f(y) \mathbb{E}_x(p_{t-\tau}(\ell,y) \chi_{\{\tau < t, W_\tau = \ell\}}) \, \mu(dy)$$

$$- \int_I f(y) \mathbb{E}_x(p_{t-\tau}(-\ell,y) \chi_{\{\tau < t, W_\tau = -\ell\}}) \, \mu(dy)$$

$$= \int_I f(y) \frac{1}{\sqrt{2\pi t}} e^{-|x-y|^2/2t} \, \mu(dy)$$

$$- \mathbb{E}_x \left( \int_I f(y) (p_{t-\tau}(\ell,y) \chi_{\{\tau < t, W_\tau = \ell\}}) \, \mu(dy) \right)$$

$$- \mathbb{E}_x \left( \int_I f(y) (p_{t-\tau}(-\ell,y) \chi_{\{\tau < t, W_\tau = -\ell\}}) \, \mu(dy) \right)$$

(where the second rearrangement follows from Fubini's theorem), and therefore

$$\int_0^t P_{t-s} f(x)\, ds$$

$$= \int_0^t \int_I f(y)\, \frac{1}{\sqrt{2\pi(t-s)}} e^{-|x-y|^2/2(t-s)}\, \mu(dy)\, ds$$

$$\quad - \mathbb{E}_x \left( \int_0^t \int_I f(y)(p_{t-s-\tau}(\ell, y) \chi_{\{\tau < t-s, W_\tau = \ell\}})\, \mu(dy)\, ds \right)$$

$$\quad - \mathbb{E}_x \left( \int_0^t \int_I f(y)(p_{t-s-\tau}(-\ell, y) \chi_{\{\tau < t-s, W_\tau = -\ell\}})\, \mu(dy)\, ds \right).$$

We can differentiate the first of these three terms with respect to $x$ under the integral sign to give

$$\int_0^t \int_I f(y) \frac{x-y}{\sqrt{2\pi}(t-s)^{3/2}} e^{-|x-y|^2/2(t-s)}\, \mu(dy)\, ds$$

$$= \int_I (x-y) f(y) \left( \int_0^t \frac{1}{\sqrt{2\pi}(t-s)^{3/2}} e^{-|x-y|^2/2(t-s)}\, ds \right) \mu(dy),$$

and now can see that this is bounded directly, since

$$(x-y) \int_0^t \frac{1}{\sqrt{2\pi}(t-s)^{3/2}} e^{-|x-y|^2/2(t-s)}\, ds$$

is bounded as $x \to y$ by the same fixed estimate as appeared in the proof of Lemma 3.1.

Bounds for the second and third integrals are proved by using just the same estimates inside the expectations $\mathbb{E}_x$, and the explicit form of the hitting time $\tau$. The result follows. $\quad\square$

The last result in this subsection gives us bounds on the spatial derivatives of solutions to $S_N$ that are *independent of $N$*, giving us a certain delicate control that will be needed later when proving convergence.

PROPOSITION 3.7.   *Suppose $v_0 \in H_0^1(I)$, $v_- \le v_0 \le v_+$ and consider* (D) *with initial conditions $(v_0, (p_{\xi,0})_{\xi \in E})$ and each $(S_N)$ with initial conditions $(v_0, \Xi_0^{(N)})$. Then there is some constant $C_4 > 0$ such that*

$$\int_0^t \|D\mathbf{v}_s\|_{L^2(I)}^2\, ds, \qquad \int_0^t \|D\mathbf{V}_s^{(N)}\|_{L^2(I)}^2\, ds \le C_4$$

*for all $N \ge 1$ and all $t \in [0, T]$.*

Proof. We write

$$\frac{d}{dt}\|\mathbf{V}_t^{(N)}\|_{L^2(I)}^2$$

$$= 2\left\langle \mathbf{V}_t^{(N)}, \frac{d}{dt}\mathbf{V}_t^{(N)}\right\rangle$$

$$= 2\langle \mathbf{V}_t^{(N)}, \Delta\mathbf{V}_t^{(N)}\rangle + \left\langle \mathbf{V}_t^{(N)}, \frac{1}{N}\sum_{i\in\mathbb{Z}\cap NI^\circ} c_{\mathbf{\Xi}_t(i)}(v_{\mathbf{\Xi}_t(i)} - \mathbf{V}_t^{(N)}(i/N))\delta_{i/N}\right\rangle$$

$$= -2\|D\mathbf{V}_t^{(N)}\|_{L^2(I)}^2 + \left\langle \mathbf{V}_t^{(N)}, \frac{1}{N}\sum_{i\in\mathbb{Z}\cap NI^\circ} c_{\mathbf{\Xi}_t(i)}(v_{\mathbf{\Xi}_t(i)} - \mathbf{V}_t^{(N)}(i/N))\delta_{i/N}\right\rangle,$$

and so, rearranging and integrating,

$$\int_0^t \|D\mathbf{V}_s^{(N)}\|_{L^2(I)}^2\,ds$$

$$= -\frac{1}{2}\int_0^t \frac{d}{ds}\|\mathbf{V}_s^{(N)}\|_{L^2(I)}^2\,ds$$

$$\quad + \frac{1}{2}\int_0^t\left\langle \mathbf{V}_s^{(N)}, \frac{1}{N}\sum_{i\in\mathbb{Z}\cap NI^\circ} c_{\mathbf{\Xi}_s(i)}(v_{\mathbf{\Xi}_s(i)} - \mathbf{V}_s^{(N)}(i/N))\delta_{i/N}\right\rangle ds$$

$$\leq \frac{1}{2}(\|\mathbf{V}_0^{(N)}\|_{L^2(I)}^2 - \|\mathbf{V}_t^{(N)}\|_{L^2(I)}^2)$$

$$\quad + \ell t\left(\max_{\xi\in E} c_\xi\right)C_3\left(C_3 + \max_{\xi\in E}|v_\xi|\right)$$

$$\leq C_3^2 + \ell T\left(\max_{\xi\in E} c_\xi\right)C_3\left(C_3 + \max_{\xi\in E}|v_\xi|\right);$$

now this right-hand side is a constant independent of $N$. The same reasoning applied to (D) gives another constant independent of $N$, and so we may take $C_4$ to be the larger of these two constants to give a simultaneous bound on

$$\int_0^t \|D\mathbf{v}_s\|_{L^2(I)}^2\,ds \quad\text{and}\quad \int_0^t \|D\mathbf{V}_s^{(N)}\|_{L^2(I)}^2\,ds. \qquad\square$$

## 4. The main convergence result.

4.1. *A decomposition.* Suppose that $(\mathbf{v},(\mathbf{p}_\xi)_{\xi\in E})$ is a solution of (D) and that $(\mathbf{V}_t,\mathbf{\Xi}_t)_{0\leq t\leq T}$ is a solution of $(S_N)$ for some $N\geq 1$ (in this subsection we largely suppress $N$ in our notation, although it will be retained later when we consider more than one value of $N$ at once).

For each $\xi \in E$ we will need to consider the process of differences $(\mathbf{C}_{\xi,N}(\boldsymbol{\Xi}_t) - \mu \llcorner \mathbf{p}_{\xi,t})_{0 \leq t \leq T}$ taking values in $H^{-1}(I)$. We will decompose this as follows:

$$\mathbf{C}_{\xi,N}(\boldsymbol{\Xi}_t) - \mu \llcorner \mathbf{p}_{\xi,t} = \mathbf{C}_{\xi,N}(\boldsymbol{\Xi}_0) - \mu \llcorner p_{\xi,0} + \int_0^t Q_{\xi,s}(\boldsymbol{\Xi}_s, \mathbf{V}_s)\, ds + M_{\xi,t},$$

where

$$Q_{\xi,s}(\Xi, V)$$
$$= \frac{1}{N} \sum_{i \in \mathbb{Z} \cap NI^\circ} \sum_{\zeta \in E \setminus \{\xi\}} (\delta_{\zeta,\Xi(i)} \alpha_{\zeta,\xi}(V(i/N)) - \delta_{\xi,\Xi(i)} \alpha_{\xi,\zeta}(V(i/N))) \delta_{i/N}$$
$$- \mu \llcorner \frac{d}{ds} \mathbf{p}_{\xi,s}$$

and the above relation is taken as the definition of the process $(M_{\xi,t})_{0 \leq t \leq T}$. We note that we have defined $Q_{\xi,s}(\Xi, V)$ for arbitrary $\Xi$, $V$, but that the functions $\mathbf{p}_{\xi,t}$ are taken as given, and so play a part in the definition.

For the purposes of this paper, we will refer to the integral of $Q_{\xi,s}(\boldsymbol{\Xi}_s, \mathbf{V}_s)$ as the *finite variation part* of this difference process and to $M_\xi$ as the *martingale part*. These names are motivated by the analogy with the definition of a Stroock–Varadhan martingale arising from a Feller process, and are justified by Lemma 4.1 below. However, as we have already remarked, here the underlying state space is not locally compact, and $M_{\xi,t}$ takes values in the space of functionals $H^{-1}$, so we need to be careful about what we mean by martingale.

LEMMA 4.1. *Suppose $\phi$ is a bounded measurable function on $I$ and consider the process*

$$(\langle \phi, \mathbf{C}_{\xi,N}(\boldsymbol{\Xi}_t) - \mu \llcorner \mathbf{p}_{\xi,t} \rangle)_{0 \leq t \leq T}.$$

*This decomposes as*

$$\langle \phi, \mathbf{C}_{\xi,N}(\boldsymbol{\Xi}_t) - \mu \llcorner \mathbf{p}_{\xi,t} \rangle = \langle \phi, \mathbf{C}_{\xi,N}(\boldsymbol{\Xi}_0) - \mu \llcorner p_{\xi,0} \rangle$$
$$+ \int_0^t \langle \phi, Q_{\xi,s}(\boldsymbol{\Xi}_s, \mathbf{V}_s) \rangle\, ds + \langle \phi, M_{\xi,t} \rangle$$

*and $(\langle \phi, M_{\xi,t} \rangle)_{0 \leq t \leq T}$ is an $(\mathcal{F}_t)_{0 \leq t \leq T}$-adapted càdlàg martingale.*

PROOF. Although this result can be made to fall under the general theory of the Stroock–Varadhan martingale (see Proposition 1.7 in Chapter 4 of Ethier and Kurtz [5]), we give the calculation here. The càdlàg property is clear. Suppose $t \in [0, T)$ and $h \in (0, T - t]$; then

$$\mathbb{E}(\langle \phi, M_{\xi,t+h} \rangle - \langle \phi, M_{\xi,t} \rangle \mid \mathcal{F}_t)$$

$$= \mathbb{E}(\langle \phi, \mathbf{C}_{\xi,N}(\mathbf{\Xi}_{t+h}) - \mathbf{C}_{\xi,N}(\mathbf{\Xi}_t) - \mu_{\llcorner}(\mathbf{p}_{\xi,t+h} - \mathbf{p}_{\xi,t}) \mid \mathcal{F}_t)$$
$$- \int_t^{t+h} \mathbb{E}(\langle \phi, Q_{\xi,s}(\mathbf{\Xi}_s, \mathbf{V}_s) \rangle \mid \mathcal{F}_t) \, ds.$$

Dividing by $h$, letting $h \downarrow 0$ and comparing with (D-prop) and ($S_N$-jump), we see at once that the derivative

$$\frac{d}{dh}\mathbb{E}(\langle \phi, M_{\xi,t+h} \rangle \mid \mathcal{F}_t)|_{h=0}$$

exists and equals 0; but now we can apply the Dominated Convergence Theorem for conditional expectation to deduce that for any $h_0 < T - t$,

$$\frac{d}{dh}\mathbb{E}(\langle \phi, M_{\xi,t+h} \rangle \mid \mathcal{F}_t)\bigg|_{h=h_0}$$
$$= \mathbb{E}\left(\frac{d}{du}\mathbb{E}(\langle \phi, M_{\xi,t+h_0+u} \rangle |\mathcal{F}_{t+h_0})\bigg|_{u=0}\bigg|\mathcal{F}_t\right) = 0,$$

and so $\mathbb{E}(\langle \phi, M_{\xi,t+h} \rangle \mid \mathcal{F}_t)$ does not depend on $h$ and therefore equals $\langle \phi, M_{\xi,t} \rangle$ for all $h$. $\quad\square$

Before leaving this subsection, we recall some of the standard machinery of jump measures and compensators in the context of our decomposition. A general treatment can be found, for example, in Chapter 22 of [14].

Given the sample paths $(\mathbf{\Xi}_t)_{0 \le t \le T}$, we define for $i \in \mathbb{Z} \cap NI^\circ$ the random *jump measures* $\kappa_i$ on $(0, T] \times E$ by

$$\kappa_i = \sum_{t \in (0,T], \, \mathbf{\Xi}_t(i) \ne \mathbf{\Xi}_{t-}(i)} \delta_{(t, \mathbf{\Xi}_t(i))},$$

and also the associated *compensators* $\nu_i$ (also measures on $(0, T] \times E$) by

$$\nu_i(dt, dy) = \sum_{\zeta \in E \setminus \{\mathbf{\Xi}_{t-}(i)\}} \alpha_{\mathbf{\Xi}_{t-}(i), \zeta}(\mathbf{V}_t(i/N)) \delta_{\zeta, y} \, dt.$$

Given these, we can now rewrite

$$Q_{\xi,s}(\Xi, V) = \frac{1}{N} \sum_{i \in \mathbb{Z} \cap NI^\circ} \left( \int_{(0,t] \times E} (\delta_{\xi,y} - \delta_{\xi, \mathbf{\Xi}_{s-}(i)}) \, \nu_i(ds, dy) \right) \delta_{i/N}$$
$$- \mu_{\llcorner}\frac{d}{ds}\mathbf{p}_{\xi,s},$$

and can express the martingale part of our decomposition as

$$M_{\xi,t} = \frac{1}{N} \sum_{i \in \mathbb{Z} \cap NI^\circ} \left( \int_{(0,t] \times E} (\delta_{\xi,y} - \delta_{\xi, \mathbf{\Xi}_{s-}(i)}) \, (\kappa_i - \nu_i)(ds, dy) \right) \delta_{i/N}$$

(note that the $\mathbf{p}_\xi$ do not enter here at all). In particular, for $\phi$ as in Lemma 4.1, we have

$$\langle \phi, M_{\xi,t} \rangle = \frac{1}{N} \sum_{i \in \mathbb{Z} \cap NI^\circ} \int_{(0,t] \times E} \phi(i/N)(\delta_{\xi,y} - \delta_{\xi,\Xi_{s-}(i)})\,(\kappa_i - \nu_i)(ds, dy).$$

The reason we have given these details is so that, when we next use this decomposition in Section 4.3, we can call on the following standard result:

LEMMA 4.2. *With $\phi$, $\nu_i$, $\kappa_i$ as above we can evaluate the $L^2$-norm of $\langle \phi, M_{\xi,t} \rangle$ thus:*

$$\mathbb{E}(\langle \phi, M_{\xi,t} \rangle^2) = \frac{1}{N^2} \sum_{i \in \mathbb{Z} \cap NI^\circ} \int_{(0,t] \times E} \phi(i/N)^2(\delta_{\xi,y} - \delta_{\xi,\Xi_{s-}(i)})^2\,\nu_i(ds, dy).$$

4.2. *The plan of campaign.* There is a standard theory of convergence and characterization of Markov processes. This is well developed for Feller processes with a locally compact state space, in which case it relies on convergence criteria for the generator of the corresponding Feller semigroup, but becomes much more complicated and less applicable in more general metric spaces. A thorough account of both cases can be found in Ethier and Kurtz [5]; in particular, the various more general convergence theorems are given in Section 4.8.

In this paper we use more hands-on estimates to prove our desired form of convergence; given our large state space, I do not know whether the result could be proved by verifying enough conditions to apply one of the above-mentioned more general convergence theorems. The ideas we will use are motivated by a treatment of Kurtz's theorem which deals solely with explicit bounds on norms and probabilities, developed in the first instance for the case of fluid limits of pure jump processes, as described in [4]. This is due to Darling and Norris; Kurtz's original argument can be found, for example, in Kurtz [15]. In a sense, we follow an infinite-dimensional version of the Darling–Norris argument in a function space; this is made possible by the smoothing diffusive properties of the time-evolution PDE that takes the place of the ODE in their theory.

Our plan of campaign for proving our main theorem is as follows:

1. Decide which quantities should converge to the deterministic behavior as $N \to \infty$, and (importantly) in what sense they should converge: we will be working mostly in certain function spaces and their duals, and the desired convergence will hold only in the appropriate topology. The relevant choices are explicit in the statement of Theorem 2.3: we measure the deviation in the potential difference functions $\mathbf{V}_t^{(N)}$ and $\mathbf{V}_t$ by the $H_0^1$-norm $\|\mathbf{V}_t^{(N)} - \mathbf{v}_t\|_{H_0^1(I)}$ and the deviation in the channel states by the $H^{-1}$-norm $\|\mathbf{C}_{\xi,N}(\boldsymbol{\Xi}_t^{(N)}) - \mu \llcorner \mathbf{p}_{\xi,t}\|_{H^{-1}(I)}$. These particular choices are not uncommon in the study of deterministic PDE. Their motivation is certainly partly that they capture the relevant sort of convergence—the process $\mathbf{C}_{\xi,N}(\boldsymbol{\Xi}_t^{(N)}) - \mu \llcorner \mathbf{p}_{\xi,t}$ can converge *only* in a fairly weak sense, since a linear combination of Dirac measures can be "close to" an absolutely continuous measure only in a weak sense—but it is also important that we can calculate using these particular norms very easily. We will see this in Section 4.5.

2. Having decided how to measure the deviations of the stochastic evolution from the deterministic, we will work (quite hard) to prove either absolute bounds or Gronwall-like growth conditions on those deviations by using properties of the equations (D) and ($S_N$) and of the function spaces involved. In fact we will prove such bounds for three different processes. First we prove an absolute bound (in probability) on the $H^{-1}$-norm of the martingale part of the difference processes $\mathbf{C}_{\xi,N}(\boldsymbol{\Xi}_t^{(N)}) - \mu \llcorner \mathbf{p}_{\xi,t}$.

3. Next we bound the growth of the $H^{-1}$-norm of the finite variation part of the difference processes $\mathbf{C}_{\xi,N}(\boldsymbol{\Xi}_t^{(N)}) - \mu \llcorner \mathbf{p}_{\xi,t}$.

4. Finally we bound the growth of the difference of potentials $\mathbf{V}_t^{(N)} - \mathbf{v}_t$. In fact, most of the work will go into bounding the $L^2$-norm $\|\mathbf{V}_t^{(N)} - \mathbf{v}_t\|_{L^2(I)}$, and then combining this with the bounds on the different parts of $\mathbf{C}_{\xi,N}(\boldsymbol{\Xi}_t^{(N)}) - \mu \llcorner \mathbf{p}_{\xi,t}$ to enable an application of Gronwall's lemma. Only after an analogous version of the main theorem has been proved with this weaker $L^2$-estimate on $\mathbf{V}_t^{(N)} - \mathbf{v}_t$ will we bootstrap our results to give the desired $H_0^1$-bound; this will follow from standard properties of the semigroup $(P_t)_{t \geq 0}$.

This plan will be executed fully in the subsections that follow.

### 4.3. *Bounding the martingale part.*

LEMMA 4.3. *For any $\phi \in L^\infty(I)$ we have*

$$\mathbb{E}(\langle \phi, M_{\xi,t} \rangle^2) \leq 8\ell \max_{\xi,\zeta \in E} \|\alpha_{\xi,\zeta}\|_\infty \frac{\|\phi\|_\infty^2}{N} t$$

*for all $t \in [0,T]$, $N \geq 1$.*

REMARK. It is for the proof of this result that we introduced jump measures and compensators in Section 4.1, as we find here that Lemma 4.2 makes our lives very much easier.

PROOF OF LEMMA 4.3. As in Lemma 4.2:

$$\mathbb{E}(\langle \phi, M_{\xi,t} \rangle^2) = \frac{1}{N^2} \sum_{i \in \mathbb{Z} \cap NI^\circ} \int_{(0,t] \times E} \phi(i/N)^2 (\delta_{\xi,y} - \delta_{\xi,\boldsymbol{\Xi}_{s-}(i)})^2 \, \nu_i(ds,dy).$$

Substituting our definition of $\nu_i$, this becomes

$$\mathbb{E}(\langle \phi, M_{\xi,t} \rangle^2)$$

$$= \frac{1}{N^2} \sum_{i \in \mathbb{Z} \cap N I^\circ} \sum_{\zeta \in E \setminus \{\boldsymbol{\Xi}_{t-}(i)\}} \int_{(0,t] \times E} \phi(i/N)^2 (\delta_{\xi,y} - \delta_{\xi,\boldsymbol{\Xi}_{s-}(i)})^2$$

$$\times \alpha_{\boldsymbol{\Xi}_{t-}(i),\zeta}(\mathbf{V}_t(i/N)) \delta_{\zeta,y} \, ds$$

$$= \frac{1}{N^2} \sum_{i \in \mathbb{Z} \cap N I^\circ} \sum_{\zeta \in E \setminus \{\boldsymbol{\Xi}_{t-}(i)\}} \int_{(0,t]} \phi(i/N)^2 (\delta_{\xi,\zeta} - \delta_{\xi,\boldsymbol{\Xi}_{s-}(i)})^2$$

$$\times \alpha_{\boldsymbol{\Xi}_{t-}(i),\zeta}(\mathbf{V}_t(i/N)) \, ds,$$

and now the given bound is clear by inspection. $\square$

LEMMA 4.4. *Fix $C > 0$ and suppose $T > 0$ and $\varepsilon > 0$. Then for any sufficiently large $N$, say $N \geq N_1$, we can find a subset $\Omega_1$ of $\Omega$ such that*

$$\mathbb{P}(\Omega \setminus \Omega_1) < \varepsilon$$

*and for any $\psi \in H_0^1(I)$ satisfying both $\|\psi\|_{L^2(I)} \leq C$ and $\|D\psi\|_{L^2(I)} < C$ we have*

$$\left( \sup_{0 \leq t \leq T} |\langle \psi, M_{\xi,t} \rangle| \right)^2 \leq \varepsilon$$

*on all of $\Omega_1$.*

REMARK. This lemma tells that for sufficiently large $N$ the martingales $\langle \psi, M_\xi \rangle$ for $\psi$ of the form described can all be controlled simultaneously with high probability. We will prove this by using estimates on this martingale for finitely many individual functions, and then approximating an arbitrary function by combinations of these.

This is possible because, by the previous lemma, we can control the size of $\langle \psi, M_{\xi,t} \rangle$ if we know only the uniform norm $\|\psi\|_\infty$ of $\psi$; but any bounded subset in the Sobolev space $H_0^1(I)$ is compact for the uniform topology on $C(I)$. [Indeed, as is well known, $H_0^1(I)$ embeds continuously into the space of functions on $I$ that are Hölder-$\frac{1}{2}$ continuous, and so we have the equicontinuity needed to apply the Arzelà–Ascoli theorem.] This means that, in the uniform norm, we can approximate the whole of any bounded subset of $H_0^1(I)$ with only finitely many of its members.

PROOF OF LEMMA 4.4. Let $E$ be the set of $\psi$ satisfying the stated bounds. Since $\langle (-\psi), M_{\xi,t} \rangle = -\langle \psi, M_{\xi,t} \rangle$ and $\psi \in E$ if and only if $-\psi \in E$, it suffices to prove the above with the last inequality replaced by

$$\left( \sup_{0 \leq t \leq T} \langle \psi, M_{\xi,t} \rangle \right)^2 \leq \varepsilon.$$

Since $(M_{\xi,t})_{0 \le t \le T}$ is a martingale, replacing $\varepsilon$ by $\frac{\varepsilon}{4}$ and using Doob's $L^2$ martingale inequality shows further that it actually suffices to prove the inequality

$$\langle \psi, M_{\xi,t} \rangle^2 \le \varepsilon.$$

STEP 1. Since we can bound $\mathbf{C}_{\xi,N}(\Xi)$, $\mu \llcorner \mathbf{p}_{\xi,t}$ and $Q_{\xi,s}(\Xi, V)$ in $H^{-1}(I)$ independently of $N$ and $t \in [0, T]$, we can choose some $\eta > 0$ such that whenever $\|\psi\|_\infty \le \eta$, then also $\langle \psi, M_{\xi,t} \rangle^2 \le \frac{\varepsilon}{4}$.

STEP 2. Use the above-mentioned compact embedding to choose finitely many $\phi_1, \phi_2, \ldots, \phi_k \in E$ so that any $\psi \in E$ has $\|\psi - \phi_j\|_\infty \le \eta$ for some $j \le k$. Applying Lemma 4.3 to each of the functions $\phi_i$, we can choose $N_1 \ge 1$ such that if $N \ge N_1$, then we have

$$\mathbb{E}(\langle \phi_i, M_{\xi,t} \rangle^2) < \frac{\varepsilon^2}{4k}$$

for all $t \in [0, T]$ and for $1 \le i \le k$. It follows from Chebyschev's inequality that

$$\mathbb{P}\left( \langle \phi_i, M_{\xi,t} \rangle^2 \ge \frac{\varepsilon}{4} \text{ for some } 1 \le i \le k, t \in [0, T] \right) \le \varepsilon;$$

set

$$\Omega_1 = \left\{ \langle \phi_i, M_{\xi,t} \rangle^2 < \frac{\varepsilon}{4} \text{ for all } 1 \le i \le k, t \in [0, T] \right\}.$$

STEP 3. Now let $\psi \in E$, and choose $j \le k$ so that $\|\psi - \phi_j\| \le \eta$. Since $\langle \psi - \phi_j, M_{\xi,t} \rangle^2 \le \frac{\varepsilon}{4}$, on $\Omega_1$ we must have

$$\langle \psi, M_{\xi,t} \rangle^2 \le (|\langle \psi - \phi_j, M_{\xi,t} \rangle| + |\langle \phi_j, M_{\xi,t} \rangle|)^2 \le \left( \frac{\sqrt{\varepsilon}}{2} + \frac{\sqrt{\varepsilon}}{2} \right)^2 = \varepsilon. \quad \square$$

COROLLARY 4.5. *Suppose $T > 0$ and $\delta, \varepsilon > 0$. Then for any sufficiently large $N$, say $N \ge N_1$, we can find a subset $\Omega_1$ of $\Omega$ such that*

$$\mathbb{P}(\Omega \setminus \Omega_1) < \varepsilon$$

*and*

$$\sup_{0 \le t \le T} \|M_{\xi,t}\|_{H^{-1}(I)} \le \delta$$

*on all of $\Omega_1$.*

4.4. *Bounded growth of the finite variation part.* In this subsection we will start to tie together the processes $\mathbf{C}_{\xi,N}(\boldsymbol{\Xi}_t^{(N)}) - \mu{\llcorner}\mathbf{p}_{\xi,t}$ and $\mathbf{V}_t^{(N)} - \mathbf{v}_t$ (retaining now the superscript $N$).

LEMMA 4.6.  *With the notation of the start of Section 3, there is a constant $C_5 > 0$ independent of $N$ such that for every $\xi \in E$, $\Xi \in E^{\mathbb{Z} \cap NI^\circ}$ and $V \in H_0^1(I)$*

$$\|Q_{\xi,t}(\Xi,V)\|_{H^{-1}(I)} \le C_5(1 + \|V\|_{H_0^1(I)}) \sum_{\zeta \in E} \|\mathbf{C}_{\zeta,N}(\Xi,V) - \mu{\llcorner}\mathbf{p}_{\zeta,t}\|_{H^{-1}(I)}$$

$$+ C_5\|V - \mathbf{v}_t\|_{L^2(I)}$$

*for all $t \in [0,T]$.*

PROOF.  Writing out more fully the definition of $Q_{\xi,t}(\Xi,V)$ and expanding using (D-prop) we have

$$Q_{\xi,t}(\Xi,V)$$

$$= \frac{1}{N} \sum_{i \in \mathbb{Z} \cap NI^\circ} \sum_{\zeta \in E \setminus \{\xi\}} (\delta_{\zeta,\Xi(i)} \alpha_{\zeta,\xi}(V(i/N)) - \delta_{\xi,\Xi(i)} \alpha_{\xi,\zeta}(V(i/N))) \delta_{i/N}$$

$$- \mu{\llcorner} \frac{d}{dt} \mathbf{p}_{\xi,t}$$

$$= \frac{1}{N} \sum_{i \in \mathbb{Z} \cap NI^\circ} \sum_{\zeta \in E \setminus \{\xi\}} (\delta_{\zeta,\Xi(i)} \alpha_{\zeta,\xi}(V(i/N)) - \delta_{\xi,\Xi(i)} \alpha_{\xi,\zeta}(V(i/N))) \delta_{i/N}$$

$$- \mu{\llcorner} \sum_{\zeta \in E \setminus \{\xi\}} ((\alpha_{\zeta,\xi} \circ \mathbf{v}_t) \cdot \mathbf{p}_{\zeta,t} - (\alpha_{\xi,\zeta} \circ \mathbf{v}_t) \cdot \mathbf{p}_{\xi,t})$$

$$= \left( \frac{1}{N} \sum_{i \in \mathbb{Z} \cap NI^\circ} \sum_{\zeta \in E \setminus \{\xi\}} \delta_{\zeta,\Xi(i)} \alpha_{\zeta,\xi}(V(i/N)) \delta_{i/N} \right.$$

$$\left. - \sum_{\zeta \in E \setminus \{\xi\}} \mu{\llcorner}(\alpha_{\zeta,\xi} \circ \mathbf{v}_t) \cdot \mathbf{p}_{\zeta,t} \right)$$

$$- \left( \frac{1}{N} \sum_{i \in \mathbb{Z} \cap NI^\circ} \sum_{\zeta \in E \setminus \{\xi\}} \delta_{\xi,\Xi(i)} \alpha_{\xi,\zeta}(V(i/N)) \delta_{i/N} \right.$$

$$\left. - \sum_{\zeta \in E \setminus \{\xi\}} \mu{\llcorner}(\alpha_{\xi,\zeta} \circ \mathbf{v}_t) \cdot \mathbf{p}_{\xi,t} \right)$$

$$= \sum_{\zeta \in E \setminus \{\xi\}} \left( \frac{1}{N} \sum_{i \in \mathbb{Z} \cap NI^\circ} \delta_{\zeta,\Xi(i)} \alpha_{\zeta,\xi}(V(i/N)) \delta_{i/N} - \mu{\llcorner}(\alpha_{\zeta,\xi} \circ \mathbf{v}_t) \cdot \mathbf{p}_{\zeta,t} \right)$$

$$- \sum_{\zeta \in E \setminus \{\xi\}} \left( \frac{1}{N} \sum_{i \in \mathbb{Z} \cap NI^\circ} \delta_{\xi, \Xi(i)} \alpha_{\xi, \zeta}(V(i/N)) \delta_{i/N} - \mu {\llcorner} (\alpha_{\zeta, \xi} \circ \mathbf{v}_t) \cdot \mathbf{p}_{\xi, t} \right).$$

This has put $Q_{\xi, t}$ into a form with which we can work: we will now consider separately the individual terms in each sum. We will show the working for the first; the second is treated similarly.

The term in question is

$$\frac{1}{N} \sum_{i \in \mathbb{Z} \cap NI^\circ} \delta_{\zeta, \Xi(i)} \alpha_{\zeta, \xi}(V(i/N)) \delta_{i/N} - \mu {\llcorner} (\alpha_{\zeta, \xi} \circ \mathbf{v}_t) \cdot \mathbf{p}_{\zeta, t};$$

let us call this $Q_{\xi, \zeta, t}(\Xi, V)$. Suppose that $\theta \in H_0^1(I)$ with $\|\theta\|_{H_0^1(I)} \leq 1$; then we find

$$\langle \theta, Q_{\xi, \zeta, t}(\Xi, V) \rangle$$

$$= \frac{1}{N} \sum_{i \in \mathbb{Z} \cap NI^\circ} \delta_{\zeta, \Xi(i)} \alpha_{\zeta, \xi}(V(i/N)) \theta(i/N) - \int_I \theta \cdot (\alpha_{\zeta, \xi} \circ \mathbf{v}_t) \cdot \mathbf{p}_{\zeta, t} \, d\mu$$

$$= \frac{1}{N} \sum_{i \in \mathbb{Z} \cap NI^\circ} \delta_{\zeta, \Xi(i)} \alpha_{\zeta, \xi}(V(i/N)) \theta(i/N) - \int_I \theta \cdot (\alpha_{\zeta, \xi} \circ V) \cdot \mathbf{p}_{\zeta, t} \, d\mu$$

$$+ \int_I \theta \cdot ((\alpha_{\zeta, \xi} \circ V) - (\alpha_{\zeta, \xi} \circ \mathbf{v}_t)) \cdot \mathbf{p}_{\zeta, t} \, d\mu.$$

The first line above is just

$$\langle \theta \cdot (\alpha_{\zeta, \xi} \circ V), \mathbf{C}_{\zeta, N}(\Xi) - \mu {\llcorner} \mathbf{p}_{\zeta, t} \rangle,$$

and so, bounding the first and second terms separately (using Cauchy–Schwarz for the second),

$$|\langle \theta, Q_{\xi, \zeta, t}(\Xi, V) \rangle| \leq \|\theta \cdot (\alpha_{\zeta, \xi} \circ V)\|_{H_0^1(I)} \cdot \|\mathbf{C}_{\zeta, N}(\Xi) - \mu {\llcorner} \mathbf{p}_{\zeta, t}\|_{H^{-1}(I)}$$

$$+ \mathrm{Lip}(\alpha_{\zeta, \xi}) \cdot \|\theta\|_\infty \cdot \sqrt{2\ell} \|V - \mathbf{v}_t\|_{L^2(I)}.$$

Next we note that for any $\phi, \psi \in H_0^1(I)$ we have

$$D(\phi \psi) = (D\phi)\psi + \phi(D\psi), \qquad \mu\text{-a.e.;}$$

this can be seen directly, using the fact that $D\psi$ and $D\phi$ both exist as the usual limit of quotients $\mu$-a.e. It follows that $\phi \psi \in H_0^1(I)$, and that we can bound $\|\phi \psi\|_{H_0^1(I)}^2$ in the following way:

$$\int_I (\phi \psi)^2 + (D(\phi \psi))^2 \, d\mu \leq \|\phi\|_\infty^2 \|\psi\|_{L^2(I)}^2 + \int_I ((D\phi)\psi + \phi(D\psi))^2 \, d\mu$$

$$\leq \|\phi\|_\infty^2 \|\psi\|_{L^2(I)}^2$$

$$+ 2 \int_I ((D\phi)\psi)^2 \, d\mu + 2 \int_I (\phi(D\psi))^2 \, d\mu.$$

Applying this with $\phi = \theta$ and $\psi = \alpha_{\zeta,\xi} \circ V$, we obtain

$$\|\theta \cdot (\alpha_{\zeta,\xi} \circ V)\|_{H_0^1(I)}^2$$
$$\leq \|\theta\|_\infty^2 \|\alpha_{\zeta,\xi} \circ V\|_{L^2}^2$$
$$+ 2 \int_I ((D\theta)(\alpha_{\zeta,\xi} \circ V))^2 + (\theta(D(\alpha_{\zeta,\xi} \circ V)))^2 \, d\mu$$
$$\leq \|\theta\|_\infty^2 \|\alpha_{\zeta,\xi} \circ V\|_{L^2}^2$$
$$+ 2\|\alpha_{\zeta,\xi}\|_\infty^2 \|D\theta\|_{L^2(I)}^2 + \|\theta\|_\infty^2 \|D(\alpha_{\zeta,\xi} \circ V)\|_{L^2(I)}^2.$$

Now by Poincaré's inequality the norm $\|\cdot\|_\infty$ is bounded by $\|\cdot\|_{H_0^1(I)}$ to within a multiplicative constant; since also $\alpha_{\zeta,\xi}$ is differentiable, $\mathrm{Lip}(\alpha_{\zeta,\xi}) < \infty$ and $\|\alpha_{\zeta,\xi}\|_\infty < \infty$, it follows that there is some $C < \infty$ for which

$$\|\theta \cdot (\alpha_{\zeta,\xi} \circ V)\|_{H_0^1(I)}^2 \leq C^2(1 + \|V\|_{L^2(I)}^2 + \|DV\|_{L^2(I)}^2) \leq C^2(1 + \|V\|_{H_0^1(I)})^2$$

whenever $\|\theta\|_{H_0^1(I)} \leq 1$. Replacing $C$ by $C \vee (\sqrt{2\ell} \max_{\zeta,\xi \in E} \mathrm{Lip}(\alpha_{\zeta,\xi}))$ if necessary and substituting back into our bound for

$$|\langle \theta, Q_{\xi,\zeta,t}(\Xi, V)\rangle|,$$

we obtain

$$|\langle \theta, Q_{\xi,\zeta,t}(\Xi, V)\rangle|$$
$$\leq C(1 + \|V\|_{H_0^1(I)})\|\mathbf{C}_{\zeta,N}(\Xi) - \mu \llcorner \mathbf{p}_{\zeta,t}\|_{H^{-1}(I)} + C\|V - \mathbf{v}_t\|_{L^2(I)}$$

when $\|\theta\|_{H_0^1(I)} \leq 1$. Summing over $\zeta \in E$ to recover the terms of our original expression for $Q_{\xi,s}(\Xi, V)$ and picking $C_5 = 2|E|C$ now gives the result. $\quad\square$

COROLLARY 4.7.    *There is a constant $C_6 > 0$ independent of $N$ such that the process $(\mathbf{C}_{\xi,N}(\mathbf{\Xi}_t^{(N)}) - \mu \llcorner \mathbf{p}_{\xi,t})_{0 \leq t \leq T}$ satisfies*

$$\|\mathbf{C}_{\xi,N}(\mathbf{\Xi}_t^{(N)}) - \mu \llcorner \mathbf{p}_{\xi,t}\|_{H^{-1}(I)}$$
$$\leq \|\mathbf{C}_{\xi,N}(\Xi_0) - \mu \llcorner p_{\xi,0}\|_{H^{-1}(I)}$$
$$+ C_6 \sqrt{\left(\int_0^t \sum_{\zeta \in E} \|\mathbf{C}_{\zeta,N}(\mathbf{\Xi}_s^{(N)}) - \mu \llcorner \mathbf{p}_{\zeta,s}\|_{H^{-1}(I)}^2 \, ds\right)}$$
$$+ C_6 \sqrt{t} \sqrt{\int_0^t \|\mathbf{V}_s^{(N)} - \mathbf{v}_s\|_{L^2(I)}^2 \, ds}$$
$$+ \|M_{\xi,t}\|_{H^{-1}(I)}.$$

PROOF.   Integrating the inequality from the previous lemma in the decomposition of the difference process yields

$$\|\mathbf{C}_{\xi,N}(\mathbf{\Xi}_t^{(N)}) - \mu_{\llcorner}\mathbf{p}_{\xi,t}\|_{H^{-1}(I)}$$

$$\leq \|\mathbf{C}_{\xi,N}(\Xi_0) - \mu_{\llcorner}p_{\xi,0}\|_{H^{-1}(I)}$$

$$+ C_5 \int_0^t \sum_{\zeta \in E} \|\mathbf{C}_{\zeta,N}(\mathbf{\Xi}_s^{(N)}) - \mu_{\llcorner}\mathbf{p}_{\zeta,s}\|_{H^{-1}(I)} \, ds$$

$$+ C_5 \int_0^t \|\mathbf{V}_s^{(N)}\|_{H_0^1(I)} \sum_{\zeta \in E} \|\mathbf{C}_{\zeta,N}(\mathbf{\Xi}_s^{(N)}) - \mu_{\llcorner}\mathbf{p}_{\zeta,s}\|_{H^{-1}(I)} \, ds$$

$$+ C_5 \int_0^t \|\mathbf{V}_s^{(N)} - \mathbf{v}_s\|_{L^2(I)} \, ds + \|M_{\xi,t}\|_{H^{-1}(I)}.$$

Now we apply the Cauchy–Schwarz inequality to each of the three integrals on the right-hand side to obtain

$$\|\mathbf{C}_{\xi,N}(\mathbf{\Xi}_t^{(N)}) - \mu_{\llcorner}\mathbf{p}_{\xi,t}\|_{H^{-1}(I)}$$

$$\leq \|\mathbf{C}_{\xi,N}(\Xi_0) - \mu_{\llcorner}p_{\xi,0}\|_{H^{-1}(I)}$$

$$+ C_5 \sqrt{t} \sqrt{\int_0^t \left( \sum_{\zeta \in E} \|\mathbf{C}_{\zeta,N}(\mathbf{\Xi}_s^{(N)}) - \mu_{\llcorner}\mathbf{p}_{\zeta,s}\|_{H^{-1}(I)} \right)^2 ds}$$

$$+ C_5 \sqrt{\int_0^t \|\mathbf{V}_s^{(N)}\|_{H_0^1(I)}^2 \, ds} \sqrt{\int_0^t \left( \sum_{\zeta \in E} \|\mathbf{C}_{\zeta,N}(\mathbf{\Xi}_s^{(N)}) - \mu_{\llcorner}\mathbf{p}_{\zeta,s}\|_{H^{-1}(I)} \right)^2 ds}$$

$$+ C_5 \sqrt{t} \sqrt{\int_0^t \|\mathbf{V}_s^{(N)} - \mathbf{v}_s\|_{L^2(I)}^2 \, ds}$$

$$+ \|M_{\xi,t}\|_{H^{-1}(I)}.$$

Another application of the Cauchy–Schwarz inequality, this time to the sum inside the first and second integrals, gives

$$\|\mathbf{C}_{\xi,N}(\mathbf{\Xi}_t^{(N)}) - \mu_{\llcorner}\mathbf{p}_{\xi,t}\|_{H^{-1}(I)}$$

$$\leq \|\mathbf{C}_{\xi,N}(\Xi_0) - \mu_{\llcorner}p_{\xi,0}\|_{H^{-1}(I)}$$

$$+ C_5 \sqrt{t} \sqrt{\int_0^t |E| \sum_{\zeta \in E} \|\mathbf{C}_{\zeta,N}(\mathbf{\Xi}_s^{(N)}) - \mu_{\llcorner}\mathbf{p}_{\zeta,s}\|_{H^{-1}(I)}^2 \, ds}$$

$$+ C_5 \sqrt{\int_0^t \|\mathbf{V}_s^{(N)}\|_{H_0^1(I)}^2 \, ds} \sqrt{\int_0^t |E| \sum_{\zeta \in E} \|\mathbf{C}_{\zeta,N}(\mathbf{\Xi}_s^{(N)}) - \mu_{\llcorner}\mathbf{p}_{\zeta,s}\|_{H^{-1}(I)}^2 \, ds}$$

$$+ C_5 \sqrt{t} \sqrt{\int_0^t \|\mathbf{V}_s^{(N)} - \mathbf{v}_s\|_{L^2(I)}^2 \, ds}$$

$$+ \|M_{\xi,t}\|_{H^{-1}(I)}.$$

This gives the desired result with

$$C_6 = C_5(1 \vee (\sqrt{T} + \sqrt{|E|C_4})),$$

where $C_4$ is the constant from Proposition 3.7 such that

$$\int_0^t \|\mathbf{V}_s^{(N)}\|_{H_0^1(I)}^2 \, ds \leq C_4$$

for all $N \geq 1$ and $t \in [0, T]$.   $\square$

The above inequality is not yet of the particular form needed to apply Gronwall's lemma (since we need a Gronwall-like bound on the growth of the *square* of the $H^{-1}$-norm of the difference process); however, this requires only one further (slightly brutal) manipulation.

COROLLARY 4.8.   *With* $C_6$ *as in the previous lemma the process* $(\mathbf{C}_{\xi,N}(\boldsymbol{\Xi}_t^{(N)}) - \mu \llcorner \mathbf{p}_{\xi,t})_{0 \leq t \leq T}$ *satisfies*

$$\|\mathbf{C}_{\xi,N}(\boldsymbol{\Xi}_t^{(N)}) - \mu \llcorner \mathbf{p}_{\xi,t}\|_{H^{-1}(I)}^2$$

$$\leq 4 \Bigg( \|\mathbf{C}_{\xi,N}(\Xi_0) - \mu \llcorner p_{\xi,0}\|_{H^{-1}(I)}^2$$

$$+ C_6^2 \int_0^t \sum_{\zeta \in E} \|\mathbf{C}_{\zeta,N}(\boldsymbol{\Xi}_s^{(N)}) - \mu \llcorner \mathbf{p}_{\zeta,s}\|_{H^{-1}(I)}^2 \, ds$$

$$+ C_6^2 t \int_0^t \|\mathbf{V}_s^{(N)} - \mathbf{v}_s\|_{L^2(I)}^2 \, ds + \|M_{\xi,t}\|_{H^{-1}(I)}^2 \Bigg).$$

PROOF.   This follows from squaring the inequality from the previous lemma and applying the Cauchy–Schwarz inequality.   $\square$

4.5.   *The full result.*   We are now able to prove the full result (Theorem 2.3):

THEOREM 4.9.   *Let* $\varepsilon > 0$, *and suppose given initial conditions* $v_0$ *and* $p_{\xi,0}$. *Then for any* $N$ *sufficiently large, say* $N \geq N_1$, *there exists an initial* *condition $\Xi_0$ for $(S_N)$ so that there is some "high-probability" $\Omega_1 \subseteq \Omega$ with $\mathbb{P}(\Omega \setminus \Omega_1) < \varepsilon$ and such that*

$$\sup_{0 \leq t \leq T} \|\mathbf{V}_t^{(N)} - \mathbf{v}_t\|_{H_0^1(I)} < \varepsilon,$$

$$\sup_{0 \leq t \leq T} \|\mathbf{C}_{\xi,N}(\mathbf{\Xi}_t^{(N)}) - \mathbf{p}_{\xi,t}\|_{H^{-1}(I)} < \varepsilon,$$

*on $\Omega_1$.*

The proof will rely on the various estimates we have made so far in the paper, and so we first recall those of the relevant constants that we will need again explicitly. There are two of these:

- $C_3$ is a uniform bound on $\sup_{0 \leq t \leq T} \|\mathbf{v}_t\|_\infty$ and $\sup_{0 \leq t \leq T} \|\mathbf{V}_t^{(N)}\|_\infty$, independent of $N$;
- $C_6$ is such that the process $(\mathbf{C}_{\xi,N}(\mathbf{\Xi}_t^{(N)}) - \mu \llcorner \mathbf{p}_{\xi,t})_{0 \leq t \leq T}$ satisfies

$$\|\mathbf{C}_{\xi,N}(\mathbf{\Xi}_t^{(N)}) - \mu \llcorner \mathbf{p}_{\xi,t}\|_{H^{-1}(I)}^2$$
$$\leq 4 \Bigg( \|\mathbf{C}_{\xi,N}(\Xi_0) - \mu \llcorner p_{\xi,0}\|_{H^{-1}(I)}^2$$
$$+ C_6^2 \int_0^t \sum_{\zeta \in E} \|\mathbf{C}_{\zeta,N}(\mathbf{\Xi}_s^{(N)}) - \mu \llcorner \mathbf{p}_{\zeta,s}\|_{H^{-1}(I)}^2 \, ds$$
$$+ C_6^2 t \int_0^t \|\mathbf{V}_s^{(N)} - \mathbf{v}_s\|_{L^2(I)}^2 \, ds + \|M_{\xi,t}\|_{H^{-1}(I)}^2 \Bigg).$$

PROOF OF THEOREM 2.3. The proof consists of a further sequence of estimates; we break it into five steps.

For the first four steps we lower our sights slightly to showing that for $N$ sufficiently large we can ensure the bounds

$$\sup_{0 \leq t \leq T} \|\mathbf{V}_t^{(N)} - \mathbf{v}_t\|_{L^2(I)} < \varepsilon,$$

$$\sup_{0 \leq t \leq T} \|\mathbf{C}_{\xi,N}(\mathbf{\Xi}_t^{(N)}) - \mathbf{p}_{\xi,t}\|_{H^{-1}(I)} < \varepsilon,$$

on some $\Omega_1$ with $\mathbb{P}(\Omega \setminus \Omega_1) < \varepsilon$; that is, our first estimate is now for $\| \cdot \|_{L^2(I)}$ rather than $\| \cdot \|_{H_0^1(I)}$. In Step 5 we will then bootstrap from this weakened result to the full theorem, using properties of the norms in question and the Feller semigroup $(P_t)_{t \geq 0}$.

STEP 1.  Calculating $\frac{d}{dt}\|\mathbf{V}_t^{(N)} - \mathbf{v}_t\|_{L^2(I)}^2$ from $(\mathrm{S}_N\text{-PDE})$ and $(\mathrm{D}\text{-PDE})$, adding the diffusion term to both sides and rearranging gives

$$\frac{d}{dt}\|\mathbf{V}_t^{(N)} - \mathbf{v}_t\|_{L^2(I)}^2 + 2\|D(\mathbf{V}_t^{(N)} - \mathbf{v}_t)\|_{L^2(I)}^2$$

$$= 2\left\langle \mathbf{V}_t^{(N)} - \mathbf{v}_t, \sum_{\xi \in E} c_\xi(\mathbf{C}_{\xi,N}(\mathbf{\Xi}_t^{(N)})) \llcorner (v_\xi - \mathbf{V}_t^{(N)}) - \mu\llcorner\mathbf{p}_{\xi,t} \cdot (v_\xi - \mathbf{v}_t)) \right\rangle$$

$$= -2\sum_{\xi \in E} c_\xi\langle \mathbf{V}_t^{(N)} - \mathbf{v}_t, \mathbf{C}_{\xi,N}(\mathbf{\Xi}_t^{(N)}) \llcorner \mathbf{V}_t^{(N)} - \mu\llcorner\mathbf{p}_{\xi,t} \cdot \mathbf{v}_t\rangle$$

$$+ 2\sum_{\xi \in E} c_\xi v_\xi \langle \mathbf{V}_t^{(N)} - \mathbf{v}_t, \mathbf{C}_{\xi,N}(\mathbf{\Xi}_t^{(N)}) - \mu\llcorner\mathbf{p}_{\xi,t}\rangle.$$

We obtain a bound on this last expression by treating the terms in these two sums separately. First we have

$$\langle \mathbf{V}_t^{(N)} - \mathbf{v}_t, \mathbf{C}_{\xi,N}(\mathbf{\Xi}_t^{(N)}) \llcorner \mathbf{V}_t^{(N)} - \mu\llcorner\mathbf{p}_{\xi,t} \cdot \mathbf{v}_t\rangle$$

$$= \langle \mathbf{V}_t^{(N)} - \mathbf{v}_t, (\mathbf{C}_{\xi,N}(\mathbf{\Xi}_t^{(N)}) - \mu\llcorner\mathbf{p}_{\xi,t}) \llcorner \mathbf{V}_t^{(N)}\rangle$$

$$+ \langle \mathbf{V}_t^{(N)} - \mathbf{v}_t, \mu\llcorner\mathbf{p}_{\xi,t} \cdot (\mathbf{V}_t^{(N)} - \mathbf{v}_t)\rangle$$

$$= \langle \mathbf{V}_t^{(N)}(\mathbf{V}_t^{(N)} - \mathbf{v}_t), \mathbf{C}_{\xi,N}(\mathbf{\Xi}_t^{(N)}) - \mu\llcorner\mathbf{p}_{\xi,t}\rangle + \langle (\mathbf{V}_t^{(N)} - \mathbf{v}_t)^2, \mu\llcorner\mathbf{p}_{\xi,t}\rangle.$$

We now bound *these* two subterms separately. The second can be bounded directly by $\|\mathbf{V}_t^{(N)} - \mathbf{v}_t\|_{L^2(I)}^2$, and the first by

$$\|\mathbf{V}_t^{(N)}(\mathbf{V}_t^{(N)} - \mathbf{v}_t)\|_{H_0^1(I)}\|\mathbf{C}_{\xi,N}(\mathbf{\Xi}_t^{(N)}) - \mu\llcorner\mathbf{p}_{\xi,t}\|_{H^{-1}(I)}$$

$$= \|(\mathbf{V}_t^{(N)} - \mathbf{v}_t)^2 + \mathbf{v}_t(\mathbf{V}_t^{(N)} - \mathbf{v}_t)\|_{H_0^1(I)}\|\mathbf{C}_{\xi,N}(\mathbf{\Xi}_t^{(N)}) - \mu\llcorner\mathbf{p}_{\xi,t}\|_{H^{-1}(I)}$$

$$\leq (\|(\mathbf{V}_t^{(N)} - \mathbf{v}_t)^2\|_{H_0^1(I)}$$

$$+ \|\mathbf{v}_t(\mathbf{V}_t^{(N)} - \mathbf{v}_t)\|_{H_0^1(I)})\|\mathbf{C}_{\xi,N}(\mathbf{\Xi}_t^{(N)}) - \mu\llcorner\mathbf{p}_{\xi,t}\|_{H^{-1}(I)}.$$

Hence we obtain for the sum of the two subterms

$$|\langle \mathbf{V}_t^{(N)} - \mathbf{v}_t, \mathbf{C}_{\xi,N}(\mathbf{\Xi}_t^{(N)}) \llcorner \mathbf{V}_t^{(N)} - \mu\llcorner\mathbf{p}_{\xi,t} \cdot \mathbf{v}_t\rangle|$$

$$\leq (\|(\mathbf{V}_t^{(N)} - \mathbf{v}_t)^2\|_{H_0^1(I)}$$

$$+ \|\mathbf{v}_t(\mathbf{V}_t^{(N)} - \mathbf{v}_t)\|_{H_0^1(I)})\|\mathbf{C}_{\xi,N}(\mathbf{\Xi}_t^{(N)}) - \mu\llcorner\mathbf{p}_{\xi,t}\|_{H^{-1}(I)}$$

$$+ \|\mathbf{V}_t^{(N)} - \mathbf{v}_t\|_{L^2(I)}^2.$$

Similarly but more straightforwardly, we have for the terms in the second sum

$$|v_\xi \langle \mathbf{V}_t^{(N)} - \mathbf{v}_t, \mathbf{C}_{\xi,N}(\boldsymbol{\Xi}_t^{(N)}) - \mu_{\llcorner} \mathbf{p}_{\xi,t} \rangle|$$
$$\leq |v_\xi| \|\mathbf{V}_t^{(N)} - \mathbf{v}_t\|_{H_0^1(I)} \|\mathbf{C}_{\xi,N}(\boldsymbol{\Xi}_t^{(N)}) - \mu_{\llcorner} \mathbf{p}_{\xi,t}\|_{H^{-1}(I)}.$$

Adding these two inequalities in our original equation and integrating with respect to $t$ gives

$$\|\mathbf{V}_t^{(N)} - \mathbf{v}_t\|_{L^2(I)}^2 + 2 \int_0^t \|D(\mathbf{V}_t^{(N)} - \mathbf{v}_t)\|_{L^2(I)}^2 \, ds$$
$$\leq \|\mathbf{V}_0^{(N)} - \mathbf{v}_0\|_{L^2(I)}^2$$
$$+ 2 \Big( \max_{\xi \in E} c_\xi \Big) \sum_{\xi \in E} \int_0^t \|(\mathbf{V}_s^{(N)} - \mathbf{v}_s)^2\|_{H_0^1(I)}$$
$$\times \|\mathbf{C}_{\xi,N}(\boldsymbol{\Xi}_s^{(N)}) - \mu_{\llcorner} \mathbf{p}_{\xi,s}\|_{H^{-1}(I)} \, ds$$
$$+ 2 \Big( \max_{\xi \in E} c_\xi \Big) \sum_{\xi \in E} \int_0^t \|\mathbf{v}_s(\mathbf{V}_s^{(N)} - \mathbf{v}_s)\|_{H_0^1(I)}$$
$$\times \|\mathbf{C}_{\xi,N}(\boldsymbol{\Xi}_s^{(N)}) - \mu_{\llcorner} \mathbf{p}_{\xi,s}\|_{H^{-1}(I)} \, ds$$
$$+ 2 \Big( \max_{\xi \in E} c_\xi \Big) \sum_{\xi \in E} \int_0^t |v_\xi| \|\mathbf{V}_s^{(N)} - \mathbf{v}_s\|_{H_0^1(I)}$$
$$\times \|\mathbf{C}_{\xi,N}(\boldsymbol{\Xi}_s^{(N)}) - \mu_{\llcorner} \mathbf{p}_{\xi,s}\|_{H^{-1}(I)} \, ds$$
$$+ 2 \Big( \max_{\xi \in E} c_\xi \Big) |E| \int_0^t \|\mathbf{V}_s^{(N)} - \mathbf{v}_s\|_{L^2(I)}^2 \, ds.$$

This is not yet in a very useful form, but to go further we will need to look first at some parts of this expression in more detail.

STEP 2. In addition to $\|\mathbf{V}_s^{(N)} - \mathbf{v}_s\|_{H_0^1(I)}$, the expressions $\|(\mathbf{V}_s^{(N)} - \mathbf{v}_s)^2\|_{H_0^1(I)}$ and $\|\mathbf{v}_s(\mathbf{V}_s^{(N)} - \mathbf{v}_s)\|_{H_0^1(I)}$ have crept into our working. We do not know any bounds on these quantities in particular, and would like to remove them altogether; it turns out that we can do this, using the fact that $\mathbf{V}_t^{(N)}$, $\mathbf{v}_t$ are uniformly bounded (Proposition 3.5). As in Section 4.4 we will use the inequality

$$\int_I (\phi\psi)^2 + (D(\phi\psi))^2 \, d\mu \leq \|\phi\|_\infty^2 \|\psi\|_{L^2}^2 + 2 \int_I ((D\phi)\psi)^2 \, d\mu + 2 \int_I (\phi(D\psi))^2 \, d\mu$$

for $\phi, \psi \in H_0^1(I)$.

We will apply this twice. In the first case we take $\psi = \phi = \mathbf{V}_t^{(N)} - \mathbf{v}_t$ to find that

$$\|(\mathbf{V}_t^{(N)} - \mathbf{v}_t)^2\|_{H_0^1(I)}^2 \leq \|(\mathbf{V}_t^{(N)} - \mathbf{v}_t)^2\|_{L^2(I)}^2$$
$$+ 4\|\mathbf{V}_t^{(N)} - \mathbf{v}_t\|_\infty^2 \|D(\mathbf{V}_t^{(N)} - \mathbf{v}_t)\|_{L^2(I)}^2$$
$$\leq \|\mathbf{V}_t^{(N)} - \mathbf{v}_t\|_\infty^2 (\|\mathbf{V}_t^{(N)} - \mathbf{v}_t\|_{L^2(I)}^2$$
$$+ 4\|D(\mathbf{V}_t^{(N)} - \mathbf{v}_t)\|_{L^2(I)}^2).$$

Second, we keep $\phi = \mathbf{V}_t^{(N)} - \mathbf{v}_t$ but take $\psi = \mathbf{v}_t$, to see that

$$\|\mathbf{v}_t(\mathbf{V}_t^{(N)} - \mathbf{v}_t)\|_{H_0^1(I)}^2$$
$$\leq \|\mathbf{v}_t(\mathbf{V}_t^{(N)} - \mathbf{v}_t)\|_{L^2(I)}^2 + 2\|D\mathbf{v}_t\|_\infty^2 \|\mathbf{V}_t^{(N)} - \mathbf{v}_t\|_{L^2(I)}^2$$
$$+ 2\|\mathbf{v}_t\|_\infty^2 \|D(\mathbf{V}_t^{(N)} - \mathbf{v}_t)\|_{L^2(I)}^2$$
$$\leq (\|\mathbf{v}_t\|_\infty^2 + 2\|D\mathbf{v}_t\|_\infty^2)\|\mathbf{V}_t^{(N)} - \mathbf{v}_t\|_{L^2(I)}^2$$
$$+ 2\|\mathbf{v}_t\|_\infty^2 \|D(\mathbf{V}_t^{(N)} - \mathbf{v}_t)\|_{L^2(I)}^2.$$

Here is the appearance of $\|D\mathbf{v}_t\|_\infty$ that we will need to bound uniformly (recall Lemma 3.6).

STEP 3. Now we can return to our estimate from the end of Step 1. Applying the Cauchy–Schwarz inequality there a few times, we obtain

$$\|\mathbf{V}_t^{(N)} - \mathbf{v}_t\|_{L^2(I)}^2 + 2\int_0^t \|D(\mathbf{V}_t^{(N)} - \mathbf{v}_t)\|_{L^2(I)}^2\,ds$$
$$\leq \|\mathbf{V}_0^{(N)} - \mathbf{v}_0\|_{L^2(I)}^2$$
$$+ 2\Big(\max_{\xi \in E} c_\xi\Big) \sum_{\xi \in E} \sqrt{\int_0^t \|(\mathbf{V}_s^{(N)} - \mathbf{v}_s)^2\|_{H_0^1(I)}^2\,ds}$$
$$\times \sqrt{\int_0^t \|\mathbf{C}_{\xi,N}(\boldsymbol{\Xi}_s^{(N)}) - \mu \llcorner \mathbf{p}_{\xi,s}\|_{H^{-1}(I)}^2\,ds}$$
$$+ 2\Big(\max_{\xi \in E} c_\xi\Big) \sum_{\xi \in E} \sqrt{\int_0^t \|\mathbf{v}_s(\mathbf{V}_s^{(N)} - \mathbf{v}_s)\|_{H_0^1(I)}^2\,ds}$$
$$\times \sqrt{\int_0^t \|\mathbf{C}_{\xi,N}(\boldsymbol{\Xi}_s^{(N)}) - \mu \llcorner \mathbf{p}_{\xi,s}\|_{H^{-1}(I)}^2\,ds}$$

$$+ 2\left(\max_{\xi \in E} c_\xi\right)\left(\max_{\xi \in E} |v_\xi|\right) \sum_{\xi \in E} \sqrt{\int_0^t \|\mathbf{V}_s^{(N)} - \mathbf{v}_s\|_{H_0^1(I)}^2 \, ds}$$

$$\times \sqrt{\int_0^t \|\mathbf{C}_{\xi,N}(\boldsymbol{\Xi}_s^{(N)}) - \mu_\llcorner \mathbf{p}_{\xi,s}\|_{H^{-1}(I)}^2 \, ds}$$

$$+ 2\left(\max_{\xi \in E} c_\xi\right)|E|\int_0^t \|\mathbf{V}_s^{(N)} - \mathbf{v}_s\|_{L^2(I)}^2 \, ds.$$

Next we will use the AM-GM inequality on the terms comprising products of square roots; however, we will apply it with a clever dodge well known in the study of PDE, observing that for any $a, b > 0$ and any $\eta > 0$ this inequality gives us

$$ab \leq \eta a^2 + \frac{1}{4\eta} b^2.$$

(See, e.g., Section 9.1 of Evans [7].)

Applying this to the above gives, for any $\eta > 0$,

$$\|\mathbf{V}_t^{(N)} - \mathbf{v}_t\|_{L^2(I)}^2 + 2\int_0^t \|D(\mathbf{V}_t^{(N)} - \mathbf{v}_t)\|_{L^2(I)}^2 \, ds$$

$$\leq \|\mathbf{V}_0^{(N)} - \mathbf{v}_0\|_{L^2(I)}^2$$

$$+ 2\left(\max_{\xi \in E} c_\xi\right)|E|\eta \int_0^t \|(\mathbf{V}_s^{(N)} - \mathbf{v}_s)^2\|_{H_0^1(I)}^2 \, ds$$

$$+ 2\left(\max_{\xi \in E} c_\xi\right) \sum_{\xi \in E} \frac{1}{4\eta} \int_0^t \|\mathbf{C}_{\xi,N}(\boldsymbol{\Xi}_s^{(N)}) - \mu_\llcorner \mathbf{p}_{\xi,s}\|_{H^{-1}(I)}^2 \, ds$$

$$+ 2\left(\max_{\xi \in E} c_\xi\right)|E|\eta \int_0^t \|\mathbf{v}_s(\mathbf{V}_s^{(N)} - \mathbf{v}_s)\|_{H_0^1(I)}^2 \, ds$$

$$+ 2\left(\max_{\xi \in E} c_\xi\right) \sum_{\xi \in E} \frac{1}{4\eta} \int_0^t \|\mathbf{C}_{\xi,N}(\boldsymbol{\Xi}_s^{(N)}) - \mu_\llcorner \mathbf{p}_{\xi,s}\|_{H^{-1}(I)}^2 \, ds$$

$$+ 2\left(\max_{\xi \in E} c_\xi\right)\left(\max_{\xi \in E} |v_\xi|\right)|E|\eta \int_0^t \|\mathbf{V}_s^{(N)} - \mathbf{v}_s\|_{H_0^1(I)}^2 \, ds$$

$$+ 2\left(\max_{\xi \in E} c_\xi\right)\left(\max_{\xi \in E} |v_\xi|\right) \sum_{\xi \in E} \frac{1}{4\eta} \int_0^t \|\mathbf{C}_{\xi,N}(\boldsymbol{\Xi}_s^{(N)}) - \mu_\llcorner \mathbf{p}_{\xi,s}\|_{H^{-1}(I)}^2 \, ds$$

$$+ 2\left(\max_{\xi \in E} c_\xi\right)|E|\int_0^t \|\mathbf{V}_s^{(N)} - \mathbf{v}_s\|_{L^2(I)}^2 \, ds.$$

Combining this with the inequalities from the end of Step 2, rearranging, and taking the chance to slim down our notation a bit, it follows that there is some fixed $C_7 > 0$ (not depending on $\eta$) such that

$$\|\mathbf{V}_t^{(N)} - \mathbf{v}_t\|_{L^2(I)}^2 + 2\int_0^t \|D(\mathbf{V}_t^{(N)} - \mathbf{v}_t)\|_{L^2(I)}^2 \, ds$$

$$\leq \|\mathbf{V}_0^{(N)} - \mathbf{v}_0\|_{L^2(I)}^2$$

$$+ \eta C_7 \int_0^t \|D(\mathbf{V}_t^{(N)} - \mathbf{v}_t)\|_{L^2(I)}^2 \, ds$$

$$+ (1+\eta) C_7 \int_0^t \|\mathbf{V}_s^{(N)} - \mathbf{v}_s\|_{L^2(I)}^2 \, ds$$

$$+ \frac{C_7}{\eta} \sum_{\xi \in E} \int_0^t \|\mathbf{C}_{\xi,N}(\mathbf{\Xi}_s^{(N)}) - \mu {\llcorner} \mathbf{p}_{\xi,s}\|_{H^{-1}(I)}^2 \, ds.$$

Now we can choose $\eta$ so small that $2 \geq \eta C_7$, and so deduce the simplified inequality

$$\|\mathbf{V}_t^{(N)} - \mathbf{v}_t\|_{L^2(I)}^2 \leq \|\mathbf{V}_0^{(N)} - \mathbf{v}_0\|_{L^2(I)}^2$$

$$+ (1+\eta) C_7 \int_0^t \|\mathbf{V}_s^{(N)} - \mathbf{v}_s\|_{L^2(I)}^2 \, ds$$

$$+ \frac{C_7}{\eta} \sum_{\xi \in E} \int_0^t \|\mathbf{C}_{\xi,N}(\mathbf{\Xi}_s^{(N)}) - \mu {\llcorner} \mathbf{p}_{\xi,s}\|_{H^{-1}(I)}^2 \, ds;$$

this follows simply by dropping the terms involving $\|D(\mathbf{V}_t^{(N)} - \mathbf{v}_t)\|_{L^2(I)}^2$. The whole point of introducing $\eta$ was to allow us to do this; now we are left with terms that we know more about.

STEP 4. Recall next our growth inequality for the finite variation parts of the difference processes from Section 4.4:

$$\|\mathbf{C}_{\xi,N}(\mathbf{\Xi}_t^{(N)}) - \mu {\llcorner} \mathbf{p}_{\xi,t}\|_{H^{-1}(I)}^2$$

$$\leq 4 \Bigg( \|\mathbf{C}_{\xi,N}(\Xi_0) - \mu {\llcorner} p_{\xi,0}\|_{H^{-1}(I)}^2$$

$$+ C_6^2 \int_0^t \sum_{\zeta \in E} \|\mathbf{C}_{\zeta,N}(\mathbf{\Xi}_s^{(N)}) - \mu {\llcorner} \mathbf{p}_{\zeta,s}\|_{H^{-1}(I)}^2 \, ds$$

$$+ C_6^2 t \int_0^t \|\mathbf{V}_s^{(N)} - \mathbf{v}_s\|_{L^2(I)}^2 \, ds + \|M_{\xi,t}\|_{H^{-1}(I)}^2 \Bigg).$$

Summing over $E$ gives

$$\sum_{\xi \in E} \|\mathbf{C}_{\xi,N}(\mathbf{\Xi}_t^{(N)}) - \mu {\llcorner} \mathbf{p}_{\xi,t}\|_{H^{-1}(I)}^2$$

$$\leq 4\Bigg(\sum_{\xi \in E} \|\mathbf{C}_{\xi,N}(\Xi_0) - \mu \llcorner p_{\xi,0}\|_{H^{-1}(I)}^2$$

$$+ C_6^2 |E| \int_0^t \sum_{\zeta \in E} \|\mathbf{C}_{\zeta,N}(\boldsymbol{\Xi}_s^{(N)}) - \mu \llcorner \mathbf{p}_{\zeta,s}\|_{H^{-1}(I)}^2 \, ds$$

$$+ C_6^2 |E| t \int_0^t \|\mathbf{V}_s^{(N)} - \mathbf{v}_s\|_{L^2(I)}^2 \, ds + \sum_{\xi \in E} \|M_{\xi,t}\|_{H^{-1}(I)}^2 \Bigg).$$

The end is now brought close with a monster application of Gronwall's lemma. Add the growth inequality obtained at the end of Step 3 to the above to find

$$\|\mathbf{V}_t^{(N)} - \mathbf{v}_t\|_{L^2(I)}^2 + \sum_{\xi \in E} \|\mathbf{C}_{\xi,N}(\boldsymbol{\Xi}_t^{(N)}) - \mu \llcorner \mathbf{p}_{\xi,t}\|_{H^{-1}(I)}^2$$

$$\leq \|\mathbf{V}_0^{(N)} - \mathbf{v}_0\|_{L^2(I)}^2$$

$$+ 4 \sum_{\xi \in E} \|\mathbf{C}_{\xi,N}(\Xi_0) - \mu \llcorner p_{\xi,0}\|_{H^{-1}(I)}^2 + 4 \sum_{\xi \in E} \|M_{\xi,t}\|_{H^{-1}(I)}^2$$

$$+ ((1+\eta)C_7 + 4C_6^2 |E| t) \int_0^t \|\mathbf{V}_s^{(N)} - \mathbf{v}_s\|_{L^2(I)}^2 \, ds$$

$$+ \left(\frac{C_7}{\eta} + 4C_6^2 |E|\right) \int_0^t \sum_{\xi \in E} \|\mathbf{C}_{\xi,N}(\boldsymbol{\Xi}_s^{(N)}) - \mu \llcorner \mathbf{p}_{\xi,s}\|_{H^{-1}(I)}^2 \, ds.$$

Letting

$$f(t) = \|\mathbf{V}_t^{(N)} - \mathbf{v}_t\|_{L^2(I)}^2 + \sum_{\xi \in E} \|\mathbf{C}_{\xi,N}(\boldsymbol{\Xi}_t^{(N)}) - \mu \llcorner \mathbf{p}_{\xi,t}\|_{H^{-1}(I)}^2,$$

we see that our above inequality implies

$$f(t) \leq A + B \int_0^t f(s) \, ds,$$

and so (by Gronwall) $f(t) \leq A \mathrm{E}^{Bt}$, where

$$A = \|\mathbf{V}_0^{(N)} - \mathbf{v}_0\|_{L^2(I)}^2 + 4 \sum_{\xi \in E} \|\mathbf{C}_{\xi,N}(\Xi_0) - \mu \llcorner p_{\xi,0}\|_{H^{-1}(I)}^2$$

$$+ 4 \sum_{\xi \in E} \sup_{0 \leq t \leq T} \|M_{\xi,t}\|_{H^{-1}(I)}^2$$

and

$$B = ((1+\eta)C_7 + 4C_6^2 |E| T) \vee \left(\frac{C_7}{\eta} + 4C_6^2 |E|\right).$$

Here $B$ does not depend on $N$, but by choosing $N$ sufficiently large we *can* make $A$ small in probability. Indeed, $\mathbf{V}_0^{(N)} = \mathbf{v}_0 = v_0$, and given $\varepsilon_1 > 0$ we note the following:

1. from Proposition 4.5, for all sufficiently large $N$, say $N \geq N_1$, we can find a subset $\Omega_1$ of $\Omega$ such that

$$\mathbb{P}(\Omega \setminus \Omega_1) < \varepsilon$$

   and

$$4 \sup_{0 \leq t \leq T} \sum_{\xi \in E} \|M_{\xi,t}\|_{H^{-1}}^2 \leq \frac{\varepsilon_1}{2}$$

   on all of $\Omega_1$;

2. for all sufficiently large $N$, say $N \geq N_2$, we can choose $\Xi_0$ so that

$$4 \sum_{\xi \in E} \|\mathbf{C}_{\xi,N}(\Xi_0) - \mu \llcorner p_{\xi,0}\|_{H^{-1}}^2 \leq \frac{\varepsilon_1}{2}$$

(this amounts to choosing $N$ so large that we may approximate each $\mu \llcorner p_{\xi,0}$ sufficiently well with a linear combination of $\delta_{i/N}$ drawn from the same $\Xi_0$; it is routine to see that this is possible). Therefore, for this choice of initial conditions and for $N$ at least $N_1 \vee N_2$, we have for all $t \in [0, T]$

$$\|\mathbf{V}_t^{(N)} - \mathbf{v}_t\|_{L^2(I)}^2, \|\mathbf{C}_{\xi,N}(\mathbf{\Xi}_t^{(N)}) - \mu \llcorner \mathbf{p}_{\xi,t}\|_{H^{-1}(I)}^2 \leq f(t) \leq \varepsilon_1 e^{BT}$$

on the large subset $\Omega_1 \subseteq \Omega$; choosing $\varepsilon_1 \leq e^{-BT}\varepsilon$ completes the proof of the weakened estimates.

STEP 5.   Finally we seek to improve our $L^2$ convergence result for the difference of potentials $\mathbf{V}^{(N)} - \mathbf{v}$ to convergence in $H_0^1$, by proving that for $N$ sufficiently large we actually have

$$\sup_{0 \leq t \leq T} \|D(\mathbf{V}_t^{(N)} - \mathbf{v}_t)\|_{L^2(I)} < \varepsilon$$

on $\Omega_1$. It turns out that this follows quickly from what we already know and the integral forms of (D-PDE) and ($S_N$-PDE). Substituting from these and subtracting, we have

$$D(\mathbf{V}_t^{(N)} - \mathbf{v}_t)$$
$$= D\left(\int_0^t P_{t-s}\left(\frac{1}{N}\sum_{i \in \mathbb{Z} \cap NI^\circ} c_{\mathbf{\Xi}_s}(v_{\mathbf{\Xi}_s^{(N)}} - \mathbf{V}_s^{(N)})\delta_{i/N}\right.\right.$$
$$\left.\left. - \sum_{\xi \in E} c_\xi \mathbf{p}_{\xi,s} \cdot (v_\xi - \mathbf{v}_s)\right) ds\right).$$

Writing

$$\mu_N = \frac{1}{N} \sum_{i \in \mathbb{Z} \cap NI^\circ} c_{\Xi_s}(v_{\Xi_s^{(N)}} - \mathbf{V}_s^{(N)})\delta_{i/N},$$

$$\mu_\infty = \mu_{\llcorner} \left( \sum_{\xi \in E} c_\xi \mathbf{p}_{\xi,s} \cdot (v_\xi - \mathbf{v}_s) \right)$$

and

$$\lambda_N = \mu_N - \mu_\infty$$

(all measures on $I$), our above expression for the derivative becomes

$$D \int_0^t P_{t-s}\lambda_N \, ds.$$

Therefore we need to prove that

$$\sup_{0 \le t \le T} \left\| D \int_0^t P_{t-s}\lambda_N \, ds \right\|_{L^2(I)} < \varepsilon$$

on $\Omega_1$ for $N$ sufficiently large.

We now deploy again our trick of breaking the integral into two parts, say over $(0, t - \varepsilon_2)$ and $(t - \varepsilon_2, t)$.

However, this time we use also two more standard observations for our Feller semigroup $(P_t)_{t \ge 0}$, which follow from the corresponding properties for the ordinary heat kernel: first, that for some fixed $C_8 > 0$ we have

$$|DP_u\lambda(x)| \le C_8 \frac{1}{\sqrt{u}} P_{u/2}\lambda(x)$$

for any finite *positive* measure $\lambda$ on $I$; and second, that there is some $C_9 > 0$ such that for any $u > 0$ we have

$$\int_I |P_u\lambda(x)|^2 \, \mu(dx) \le C_9 \frac{1}{\sqrt{u}} \|\lambda\|$$

[recalling that $\mu$ denotes Lebesgue measure, and writing $\|\lambda\|$ for the usual norm of a measure regarded as a linear functional on $C(I)$]. Combining these allows us to estimate the "awkward" part of our integral: for $t - \varepsilon_2 < s < t$ and any finite positive measure $\lambda$ we have

$$\|DP_{t-s}\lambda\|_{L^2(I)} = \sqrt{\int_I |DP_{t-s}\lambda(x)|^2 \, \mu(dx)}$$

$$\le C_8 C_9 \frac{1}{\sqrt{t-s}} \sqrt{\int_I |P_{(t-s)/2}\lambda(x)|^2 \, \mu(dx)}$$

$$\le 2C_8 C_9 \frac{1}{(t-s)^{3/4}} \|\lambda\|.$$

Hence, by the integral triangle inequality,

$$\left\| D \int_{t-\varepsilon_2}^{t} P_{t-s} \lambda_N \, ds \right\|_{L^2(I)} \leq 2C_8 C_9 \|\lambda_N\| \int_0^{\varepsilon_2} \frac{1}{u^{3/4}} \, du = 8C_8 C_9 \|\lambda_N\| \varepsilon_2^{1/4}.$$

Since $\|\lambda_N\| \leq \|\mu_N\| + \|\mu_\infty\|$ is bounded uniformly in $N$ (since $\mathbf{V}_t^{(N)}$, $\mathbf{v}_t$ are bounded uniformly in $N$ for $t \in [0,T]$), we can choose $\varepsilon_2$ so small that

$$\left\| D \int_{t-\varepsilon_2}^{t} P_{t-s} \lambda_N \, ds \right\|_{L^2(I)} \leq \tfrac{1}{2}\varepsilon$$

for every $N$, every $t \in [0,T]$ and on the whole state space $\Omega$.

On the other hand, owing to the smoothing properties of $P_{t-s}$ for $s < t$ bounded away from $t$, if we choose $N$ large enough that

$$\sup_{0 \leq t \leq T} \|\mathbf{V}_t^{(N)} - \mathbf{v}_t\|_{L^2(I)}$$

and

$$\sup_{0 \leq t \leq T} \|\mathbf{C}_{\xi,N}(\mathbf{\Xi}_t^{(N)}) - \mathbf{p}_{\xi,t}\|_{H^{-1}(I)}$$

are sufficiently small on the whole of the high-probability event $\Omega_1$, then we can ensure that we also have

$$\left\| D \int_0^{t-\varepsilon_2} P_{t-s} \lambda_N \, ds \right\|_{L^2(I)} \leq \tfrac{1}{2}\varepsilon$$

for any $t \in [0,T]$ and on the whole of $\Omega_1$. Combining these two estimates now gives the final result.  $\square$

## 5. Closing remarks.

5.1. *Appropriateness of the model.*   It is worth remarking on other stochastic models of nerve axons. These have tended to concentrate on using a white noise continuously distributed along the axon to model its stochastic nature, coupled via a suitable nonlinear parabolic PDE to the potential difference $V$ in the same way as for our stochastic individual ion channels. This white-noise approach leads to a more traditional system of SPDE, for which there are corresponding existence and uniqueness results. A good introduction to this approach is given starting at Chapter 6 of Tuckwell [18], who also describes various stochastic approximations that save on computational expense.

However, it is not clear that the regime in which this SPDE model is really appropriate is very physically interesting, for in this regime the channel size must be considered negligible even though the fluctuations caused by their stochastic nature are not (this is the regime in which spatial white noise arises). It is arguably more natural to build a model that takes individual ion channels into account and then, if desired, proceed directly to the deterministic limit given by the classical Hodgkin–Huxley equations. This is the procedure usually followed in physiology, and with which the present paper is concerned.

5.2. *Some further directions.* We finish by describing some further directions in which the rigorous analysis of the stochastic Hodgkin–Huxley equations might be taken:

1. The existence and convergence results proved above are of more academic than computational interest. However, they were originally motivated by a rather more practical problem.

   There has recently been growing interest in the deviation of a real axon from the deterministic behavior predicted by the Hodgkin–Huxley theory as a result of the stochastic nature of its components. In particular, Faisal, White and Laughlin [8] have investigated numerically the question of whether a sufficiently small axon might suffer frequent spontaneous action potentials generated by the chance event that a small number of $Na^+$ channels in close proximity stay open longer than usual and so cause a small initial rise in the membrane potential in their vicinity. They find that this can occur with a probability that increases greatly as the axonal diameter drops below about $0.1 \ \mu m$.

   Faisal, White and Laughlin's approach uses a purpose-coded computer simulator of axonal behavior and has a very high computational expense. It would be valuable to have cleaner, analytic bounds on the probabilities (or, equivalently, the long-term rates) of such events occurring. One might conjecture that, in a suitable regime of many small ion channels but with time speeded up accordingly, spontaneous action potentials appear distributed roughly as a Poisson point process in the relevant space-time band $[0, T] \times I$, where $T$ is taken very large.

   Such a result could appear as a sort of large-deviation principle around a point of equilibrium for our whole system, analogous to the analysis of metastability through large-deviation theory for finite-dimensional dynamical systems perturbed by a weak additive noise, as developed by Freidlin and Wentzell in [10]; however, their techniques would need some adaptation to suit the case of an infinite-dimensional system coupled to a large discrete system, as in the stochastic Hodgkin–Huxley model. The methods in the present paper do not seem to extend so far.

   In fact, this problem of determining the rate of spontaneous action potential generation fits naturally among various other questions that can be asked about the still-stochastic behavior of a real axon. Such a development in the special case of this model might also take into account the different relative effects of the noise in the system while it is undergoing different behavior: in particular, it seems that for realistic channel-numbers on patches of the axonal membrane there is much greater noise around local sub-threshold behavior (as when at equilibrium) than super-threshold potentials (as during the transmission of the front of an action potential). Relatedly, this analysis might also consider the possible deterioration of an action potential in the stochastic model, for which a strict solitary wave solution may not exist. Finally, we mention that Steinmetz, Manwani and Koch [17] have recently studied the reliability in the times of the spikes output by a neuron considered as a transmission of the input times by modeling a small patch of ion channels using the (nonspatial) stochastic Hodgkin–Huxley model (they also conduct a similar investigation of an alternative, the "integrate-and-fire" model). As they remark more generally,

> "The stochastic Markov version of the HH model converges to the classical, deterministic model as the number of channels grows large, but for realistic channel numbers, the stochastic model can exhibit a wide variety of behaviors (spontaneous spiking, bursting, chaos, and so on) that cannot be observed in deterministic model . . . "

It remains to be seen to what extent a more analytic treatment can be given of these different behaviors.

2. The model analyzed in this paper ignores possible external effects acting on the axon. A separate task which might be of interest is to perform the analysis of convergence in the case of an axon subject to some stimulus: for example, the arrival of a signal from the soma along the axon, modeled by changing the boundary conditions of our PDE (so that a particular input arrives at the soma end, and the conditions at the other end are free), or the application of a trans-membrane current along the length of the axon, as is sometimes used in experiments to stimulate an action potential. These stimuli could themselves be deterministic or stochastic; in the former case, one would expect the behavior of the stochastic model to converge to the trajectory of an appropriately modified PDE, while in the latter, even the limit model would have stochastic components.

For deterministic inputs, it seems likely that the methods of the present paper could still be brought to bear. However, the details of the estimates needed for the proofs both of existence and regularity and then of convergence might become considerably more complicated, and we have not tried to work out the details. In the case of stochastic inputs, some further modification of the "hands-on" convergence machinery of Darling and Norris [4], as we have adapted it to our needs in this paper, would be needed, possibly using a suitable coupling of the full stochastic model and the limiting stochastic model as $N \to \infty$.

3. A different modification of the models studied in this paper would be to work over a membrane regarded as a two-dimensional surface. This could both increase the accuracy of the model studied here (recall our early simplifying assumption to treat the axon as a line segment rather than a cylinder), and possibly give an analogous account in the case of other, two-dimensional excitable membranes that appear in physiology (see, e.g., Hille [11]). As for the second extension mentioned above, it seems likely that the basic methods of this paper would still be useful, but in this case the required estimates would probably become much more difficult, and depend strongly on the underlying two-dimensional geometry (as far as we know, even the purely deterministic model has not been rigorously analyzed in this setting).

**Acknowledgments.** My thanks go to Dr. James Norris (Department of Pure Mathematics and Mathematical Statistics, University of Cambridge) for the motivation to study this problem in the first place and for several helpful discussions and suggestions afterwards. I am grateful also to Prof. Craig Evans (Department of Mathematics, University of California at Berkeley) for advice on the existence theorem for the stochastic Hodgkin–Huxley equations, and to an anonymous reviewer for suggesting several improvements to the paper.

## REFERENCES

[1] CHOW, C. C. and WHITE, J. A. (1996). Spontaneous action potentials due to channel fluctuations. *Biophys. J.* **71** 3013–3021.

[2] CRONIN, J. (1987). *Mathematical Aspects of Hodgkin–Huxley Neural Theory.* Cambridge Univ. Press. MR0909892

[3] DEFELICE, L. J. and ISAAC, A. (1992). Chaotic states in a random world. *J. Stat. Phys.* **70** 339–352.

[4] DARLING, R. W. R. Fluid limits of pure jump Markov processes: A practical guide. Available at arXiv:math.PR/0210109.

[5] ETHIER, S. N. and KURTZ, T. G. (1986). *Markov Processes*: *Characterization and Convergence.* Wiley, New York. MR0838085

[6] EVANS, J. and SHENK, N. (1970). Solutions to axon equations. *Biophys. J.* **10** 1090–1101.

[7] EVANS, L. C. (1998). *Partial Differential Equations.* Amer. Math. Soc., Providence, RI. MR1625845

[8] FAISAL, A. A., WHITE, J. A. and LAUGHLIN, S. B. (2005). Ion-channel noise places limits on the miniaturization of the brain's wiring. *Current Biology* **15** 1143–1149.

[9] FOX, R. F. and LU, Y. (1994). Emergent collective behavior in large numbers of globally coupled independently stochastic ion channels. *Physical Review E* **49** 3421–3431.

[10] FREIDLIN, M. I. and WENTZELL, A. D. (1998). *Random Perturbations of Dynamical Systems*, 2nd ed. Springer, New York. MR1652127

[11] Hille, B. (2001). *Ion Channels of Excitable Membranes.* Sinauer Associates, Sunderland.

[12] Hodgkin, A. L. and Huxley, A. F. (1952). A quantitative description of membrane current and its application to conduction and excitation in nerve. *J. Physiol.* **117** 500–544.

[13] Jacod, J. and Shiryaev, A. N. (1987). *Limit Theorems for Stochastic Processes.* Springer, Berlin. MR0959133

[14] Kallenberg, O. (2002). *Foundations of Modern Probability*, 2nd ed. Springer, New York. MR1876169

[15] Kurtz, T. G. (1981). *Approximation of Population Processes.* SIAM, Philadelphia. MR0610982

[16] Lamberti, L. (1986). Solutions to the Hodgkin–Huxley equations. *Appl. Math. Comput.* **18** 43–70. MR0815772

[17] Steinmetz, P. N., Manwani, A. and Koch C. (2001). Variability and coding efficiency of noisy neural spike encoders. *BioSystems* **62** 87–97.

[18] Tuckwell, H. C. (1989). *Stochastic Processes in the Neurosciences.* SIAM, Philadelphia. MR1002192

Department of Mathematics
University of California at Los Angeles
Los Angeles, California 90095-1555
USA
E-mail: timaustin@math.ucla.edu
URL: www.math.ucla.edu/˜timaustin/